\documentclass[a4paper, 12pt]{article}
\usepackage{multicol}
\usepackage{amssymb,amsmath,amsthm,bm}
\usepackage[colorlinks=true, urlcolor=blue,linkcolor=blue, citecolor=blue]{hyperref}
\usepackage{mathrsfs,verbatim}
\usepackage{caption,cleveref }
\usepackage{graphicx, float}
\usepackage{color, mathtools, tikz, xcolor}
\usepackage{enumerate,enumitem}
\usepackage{tikz}
\usetikzlibrary{shapes.geometric}
\usepackage{bbm}

\usetikzlibrary{calc}

\newcommand{\sG}{lexical }

\usepackage[mathscr]{euscript}
\usepackage{appendix}
\usepackage[margin=2cm]{geometry}

\usepackage{enumitem}

\setlist[enumerate,1]{label=\upshape(\roman*),ref=(\roman*)}
\usepackage{enumitem} 
\setlist[enumerate]{topsep=0pt,itemsep=-1ex,partopsep=1ex,parsep=1ex}
\setlist[itemize]{topsep=0pt,itemsep=-1ex,partopsep=1ex,parsep=1ex}

\numberwithin{equation}{section}

	\newtheorem{theorem}{Theorem}[section]
	\newtheorem{problem}[theorem]{Problem}
	
	\newtheorem{corollary}[theorem]{Corollary}

	\newtheorem{claim}[theorem]{Claim}
	
	\newtheorem{lemma}[theorem]{Lemma}

\newenvironment{proofclaim}[1][Proof of claim]{\begin{proof}[#1]}{\end{proof}}

\def\eps{{\varepsilon}}
\def\sm{{\setminus}}

\def\ar{{\mathrm{ar}}}
\def\sar{{\mathrm{ar}^{\star}}}
\def\ea{{\mathrm{ea}}}
\def\va{{\mathrm{va}}}
\def\nsar{{\mathrm{n}\sar}}

\def\lex{{\mathrm{lex}}}
\def\part{{\mathrm{par}}}
\def\cen{\circledast}
\def\wt{\widetilde}

\newcommand{\ova}[1]{\overrightarrow{#1}}

\newcommand{\mc}[1]{\mathcal{#1}}
\newcommand{\ms}[1]{\mathscr{#1}}

\title{Rainbow subgraphs of star-coloured graphs}

\author{Allan Lo\thanks{School of Mathematics, University of Birmingham, B15 2TT, United Kingdom. Email: s.a.lo@bham.ac.uk.
AL was partially supported by EPSRC, grant no. EP/V048287/1.
} \and Klas Markstr\"om\thanks{Department of Mathematics and Mathematical Statistics, Ume\r{a} Universitet, 901 87 Ume\r{a}, Sweden. Email: klas.markstrom@umu.se.}  \and Dhruv Mubayi\thanks{Department of Mathematics, Statistics, and Computer Science, University of Illinois
Chicago, IL 60607, USA. Email: mubayi@uic.edu.} \and Katherine Staden\thanks{School of Mathematics and Statistics, The Open University, Walton Hall, Milton Keynes, UK. Email: katherine.staden@open.ac.uk. KS was supported by EPSRC Fellowship EP/V025953/1.} \and Maya Stein\thanks{Center for Mathematical Modeling (CMM - IRL CNRS 2807) and Department for Mathematical Engineering, University of Chile. Email: mstein@dim.uchile.cl. MS acknowledges support by FONDECYT Regular Grant 1221905 and by CMM BASAL FB210005 for Centers of Excellence  (ANID-Chile).} \and Lea Weber\thanks{Division of Theoretical Systems Biology, German Cancer Research Center (DKFZ), Heidelberg, Germany. Email: lea.weber@dkfz.de.}}
\date{\today}

\begin{document}

\setlength{\abovedisplayskip}{3pt}
\setlength{\belowdisplayskip}{3pt}

\maketitle

\begin{abstract}
An edge-colouring of a graph $G$ can fail to be rainbow for two reasons: either it contains a monochromatic cherry (a pair of incident edges), or a monochromatic matching of size two. A colouring is a proper colouring if it forbids the first structure, and a star-colouring if it forbids the second structure. In this paper, we study rainbow subgraphs in star-coloured graphs and determine the maximum number of colours in a star-colouring of a large complete graph which does not contain a rainbow copy of a given graph $H$. This problem is a special case of one studied by Axenovich and Iverson on generalised Ramsey numbers and we extend their results in this case.
\end{abstract}

\section{Introduction}

Given a graph $G$, let $v(G) = |V(G)|$ and $e(G)=|E(G)|$. An \emph{edge-colouring} of a graph $G$ is a function which maps its edge set $E(G)$ to a set of colours. We say that an edge-coloured graph $H$ is \emph{rainbow} if each edge receives a different colour.
An edge-colouring of a graph $G$ can fail to be rainbow for two reasons: either it contains a monochromatic \emph{cherry} $K_{1,2}$ or it contains a monochromatic matching of size 2.   Edge-colourings without the first structure are called \emph{proper}. In edge-colourings which avoid monochromatic matchings, every colour class is either a star or a triangle. We will call such colourings without triangles \emph{star-colourings}. 
Before discussing these two types of colourings, we take a step back to consider the more general question of which subgraphs appear in \emph{any} edge-coloured complete graph. 
This is answered by the famous \emph{canonical Ramsey theorem} of Erd\H{o}s and Rado~\cite{ErdosRado}.

\begin{theorem}[The canonical Ramsey theorem~\cite{ErdosRado}]\label{th:canonical}
    For any positive integer $k$, any sufficiently large edge-coloured complete graph contains a clique $K_k$ which either
    \begin{itemize}
        \item is monochromatic,
        \item is rainbow,
        \item has a \emph{lexical colouring}; that is, its vertex set can be ordered as $v_1<\cdots < v_k$  such that for $i<j$ and $i'<j'$, the colour of $v_iv_j$ is the same as the colour of  $v_{i'}v_{j'}$ if and only if $i=i'$.
    \end{itemize}
\end{theorem}

Rainbow subgraphs of properly edge-coloured graphs have been extensively studied. 
Since proper colourings avoid both  
monochromatic and
lexically coloured cliques $K_k$ (unless $k\le 2$),
Theorem~\ref{th:canonical}  implies that given any graph $H$, every sufficiently large properly coloured complete graph  contains a rainbow copy of $H$.
Less is known about graphs $H$ whose size is comparable with that of the host graph.
Since a proper colouring of $K_n$ can have as few as $n-1$ colours, we can only hope to find rainbow copies of $H$ with $e(H) \leq n-1$.
Regarding long paths, Maamoun and Meyniel~\cite{MaamounMeyniel} found an example of a proper colouring of $K_n$ without a rainbow Hamilton path for $n \geq 4$ a power of $2$. In 1989, Andersen~\cite{Andersen} conjectured that one can however guarantee a rainbow path of length $n-2$.
Even finding a rainbow path of length $n-o(n)$ was open for a long time, and was eventually achieved by Alon, Pokrovskiy and Sudakov in 2017~\cite{AlonPokrovskiySudakov}.
It has been shown that one can guarantee almost spanning cycles~\cite{AlonPokrovskiySudakov}, trees~\cite{MontgomeryPokrovskiySudakov}, and bounded degree graphs~\cite{EhardGlockJoos}. See the survey~\cite{PokrovskiySurvey} of Pokrovskiy for further results.

Now, we turn to edge-colourings of $K_n$ that are not necessarily proper. If they use sufficiently many colours, one can still guarantee a rainbow copy of a given graph $H$. The \emph{anti-Ramsey number} $\ar(n,H)$ of $H$, introduced by Erd\H{o}s, Simonovits and S\'os in 1975~\cite{ErdosSimonovitsSos}, is defined as$$
\ar(n,H) := \max\{s \in \mathbb{N}: \text{there exists an }s\text{-coloured }K_n\text{ without a rainbow copy of }H\}.
$$
Define the \emph{extremal number} ${\rm ex}(n,\mc{H})$ of a family $\mc{H}$ of graphs as
$$
{\rm ex}(n,\mc{H}) := \max\{e: \text{there exists an $n$-vertex graph  with }e\text{ edges and no copy of any }H \in \mc{H}\}
$$
and let ${\rm ex}(n,H) := {\rm ex}(n,\{H\})$.
It has been shown that
\begin{gather}
%\nonumber 
\label{eq:exantiramsey1}
\ar(n,C_k) = \left(\frac{k-2}{2}+\frac{1}{k-1}\right)n+O(1),\quad
\ar(n,K_k) = {\rm ex}(n,K_{k-1})+1,
\end{gather}
%\quad\text{and}\\
where $C_k$ is the cycle of length $k$, and further,
\begin{gather}
\label{eq:exantiramsey}
{\rm ex}(n,\mc{F}^-_H)+1 \leq \ar(n,H) \leq {\rm ex}(n,H)
\end{gather}
for any graph $H$, where $\mc{F}^-_H := \{H \setminus e: e \in E(H)\}$.
The first result is due to Montellano-Ballesteros and Neumann-Lara~\cite{MontellanoNeumann}, proving a conjecture from~\cite{ErdosSimonovitsSos},
and the other results are due to~\cite{ErdosSimonovitsSos}.
Anti-Ramsey numbers have been studied for many more $H$, see the dynamic survey~\cite{RainbowSurvey} for a more comprehensive overview.

In this paper, we initiate the study of rainbow subgraphs of star-coloured graphs.
In particular we are interested in 
%the extremal problem of 
 determining the maximum number of colours in a star-colouring of the complete graph which avoids a rainbow copy of a given graph $H$.
We say that a colouring is an \emph{$s$-star-colouring} if it is a star-colouring with $s$ colours.
Given a graph $H$ and an integer~$n$, we define the \emph{star-anti-Ramsey number} $\sar(n,H)$ of $H$ as
$$
\sar(n,H) := \max\{s \in \mathbb{N}: \text{there exists an }s\text{-star-coloured }K_n\text{ without a rainbow copy of }H\}.
$$
Clearly, we have
\begin{equation}\label{eq:antiramsey}
\sar(n,H) \leq \ar(n,H).
\end{equation}
If $n < v(H)$, then $\sar(n,H)=\ar(n,H) = \binom{n}{2}$, so we only consider the case $n \geq v(H)$.
Note that the third type of colouring in Theorem~\ref{th:canonical} is a star-colouring, so this theorem does not imply bounds on $\sar(n,H)$ for general $H$.
%does not, for general $H$, tell us anything about $\sar(n,H)$.

\subsection{Background}\label{sec:background}

We are not aware of $\sar(n,H)$ being studied directly before.
Our starting point is a closely-related result of Axenovich and Iverson~\cite{AxenovichIverson} which demonstrates a connection between $\sar(n,H)$ and the so-called vertex arboricity of $H$.
Before stating it, we recall the famous theorem of Erd\H{o}s and Stone~\cite{ErdosStone} relating the extremal number ${\rm ex}(n,H)$ and the chromatic number $\chi(H)$. (Part~(ii) is a weak consequence of another famous theorem of K\H{o}v\'ari, S\'os and Tur\'an~\cite{KovariSosTuran} which we state later in Theorem~\ref{th:KST}.)
The \emph{chromatic number} $\chi(H)$ of $H$ is the minimum number of parts in a vertex partition of $H$ such that $H$ induces an %empty 
edgeless graph in each part.

\begin{theorem} [\cite{ErdosStone, KovariSosTuran}]
Let $H$ be a graph and let $n$ be sufficiently large. 
\begin{itemize}
    \item[(i)]  If $\chi(H) \geq 3$, then
    $$
    {\rm ex}(n,H) = (1+o(1))\left(1-\frac{1}{\chi(H)-1}\right)\binom{n}{2}.
    $$
    \item[(ii)]  If $\chi(H) \leq 2$, then there exists $c=c(H)>0$ such that ${\rm ex}(n,H) \leq n^{2-1/c}$.
\end{itemize}
\end{theorem}

A huge literature has developed studying ${\rm ex}(n,H)$ for bipartite $H$ (that is, with chromatic number at most $2$), and more precise bounds than what is stated above are known, as well as even more precise results for specific $H$ of interest. We will require some of these results later.

Let $\va(H)$ denote the \emph{vertex arboricity} of $H$, which is the minimum number of parts in a vertex partition of $H$ such that $H$ induces a forest in each part.
Axenovich and Iverson~\cite{AxenovichIverson} proved that, in a broad sense, $\sar(n,H)$ and $\va(H)$ are related to each other in the same way as ${\rm ex}(n,H)$ and $\chi(H)$.

\begin{theorem}[\cite{AxenovichIverson}]\label{th:AxenovichIverson}
Let $H$ be a graph and let $n$ be sufficiently large. Then the following hold.
\begin{itemize}
    \item[(i)] If $\va(H) \geq 3$, then
    $$
    \sar(n,H) = (1+o(1))\left(1-\frac{1}{\va(H)-1}\right)\binom{n}{2}.
    $$
    \item[(ii)] If $\va(H) \leq 2$, then there exists $c=c(H)>0$ such that $\sar(n,H) \leq n^{2-1/c}$.
\end{itemize}
\end{theorem}

We explain how this result follows from~\cite{AxenovichIverson}.
Let $\mc{F},\mc{H}$ be families of graphs and let $n$ be a positive integer. Define the \emph{generalised Ramsey number} $R_{\mc{F}}(n,\mc{H})$ 
(sometimes called the \emph{mixed Ramsey number}) to be the maximum number of colours in an edge-colouring of $K_n$ without a monochromatic copy of any $F \in \mc{F}$ or a rainbow copy of any $H \in \mc{H}$.
Thus
$$
\sar(n,H) = R_{\{M_2,K_3\}}(n,\{H\})
$$
where $M_2$ is the matching of size two.
Axenovich and Iverson~\cite{AxenovichIverson} showed that, whenever $F$ is not a star, the conclusion of Theorem~\ref{th:AxenovichIverson} holds with $R_{\{F\}}(n,\{H\})$ replacing $\sar(n;H)$. We have $\sar(n,H) \leq R_{\{M_2\}}(n,\{H\})$, while there is a lower bound construction for $R_{\{M_2\}}(n,\{H\})$ which is a star-colouring; thus we obtain Theorem~\ref{th:AxenovichIverson}.

In conclusion, we know the correct order of magnitude of $\sar(n,H)$ when $\va(H) \geq 3$, while very little is known when $\va(H) \leq 2$.
The purpose of this paper is to examine the regime $\va(H) \leq 2$.

Before proceeding further, it is helpful to make two simple observations (see Lemma~\ref{lm:obs1} for their short proofs).
Let $n\in\mathbb N$ and let  $H$ be a graph. 
\medskip
\begin{itemize}
    \item Every star-colouring of $K_n$ has at least $n-1$ colours.
    \item For sufficiently large $n$, $\sar(n,H)$ exists if and only if $\va(H) \geq 2$.
\end{itemize}

\medskip
Thus, we need only consider those $H$ with $\va(H)=2$ (and such $H$ satisfy $\sar(n,H) \geq n-1$).

\medskip
For more related work, see~\cite{AxenovichKundgen,CzygrinowMollaNagle,JamisonJiangLing,JungicKaiserKral,Li}.

\subsection{Results}

Now we state our main results, starting with (mainly sharp) results for specific $H$ fulfilling $\va(H)=2$.
First, note that clearly, $\va(C_k)=2$ for any integer $k \geq 3$. 
Our first main result is for cycles of all lengths. Write $\ova{C_k}$ for the oriented graph with vertex set $\{v_1, \ldots, v_k\}$ and arc set $\{\ova{v_iv_{i+1}}: i=1, \ldots, k\}$, where subscripts are taken modulo $k$.
Recall that a tournament is an orientation of a complete graph.

\begin{theorem}\label{th:cycle}
	For all integers $n \geq k \geq 3$ we have $\sar(n,C_k)=n+{k-2 \choose 2}-1$.
    
    Moreover, every rainbow $C_k$-free star-colouring of $K_n$ with $n+\binom{k-2}{2}-1 =n-k+1+{k-1\choose 2}$ colours has the following structure. Let $A,B$ be disjoint vertex sets of sizes $n-k+1,k-1$ respectively and let $\ova{T}$ be an oriented graph obtained from a $\ova{C_k}$-free tournament on $A$ and adding all possible directed edges from~$A$ to~$B$. 
    For each vertex $x \in A$ and $y \in N_{\ova{T}}^+(x)$, colour $\ova{xy}$ with $x$. Give each pair in $B$ a new distinct colour.
\end{theorem}

In particular, Theorem~\ref{th:cycle} covers the case $C_3=K_3$, giving $\sar(n,K_3)=n-1$. Since for $K_k$ with $k\ge 5$ we have that $\va(K_k) \geq 3$ (and thus Theorem~\ref{th:AxenovichIverson} applies), 
the only clique for which the approximate value of $\sar(n,K_k)$ does not follow from existing results is $K_4$. We give a sharp result for this case.

\begin{theorem}\label{th:K4}
For all $n \ge 4$, $\sar(n, K_4) = 2n-3$.
\end{theorem}
This theorem (together with~\eqref{eq:exantiramsey1}) shows that $\ar(n,H)$ and $\sar(n,H)$ can be very far apart. 
There are many extremal colourings for $K_4$, and we are unable to characterise them, but describe some such colourings in Section~\ref{sec:K4lower}. This perhaps explains why the proof of Theorem~\ref{th:K4} is by far the most involved of the paper.

Given an integer $k \geq 3$, denote by $K_k^-$ the graph obtained from $K_k$ by removing one edge.
We  prove sharp bounds for $\sar(n,K_4^-)$ and $\sar(n,K_5^-)$. Note that both $K_4^-$ and $K_5^-$ have vertex arboricity $2$.

\begin{theorem}\label{th:K4-}
For all $n \geq 2$, $\sar(n, K_4^-) = \lfloor 3(n-1)/2 \rfloor$.

Moreover, every rainbow $K_4^-$-free star-colouring of $K_n$ with $\lfloor 3(n-1)/2\rfloor$ colours has the following structure.
There is an ordering $v_1,\ldots,v_n$ of its vertices such that we colour $v_iv_j$ with colour $\min\{i,j\}$,
and
then we give a new, distinct colour to each component in the 
following graph $J$:
\begin{itemize}
    \item if $n=2m+1$ is odd, then $E(J) = \{ v_1v_2, v_3v_4,\ldots,v_{2m-1}v_{2m}\}$;
    \item if $n=2m$ is even, then $$E(J) =  \{v_1v_2,\ldots,v_{2a-3}v_{2a-2} \} \cup E(S) \cup \{v_{2a+2}v_{2a+3},\ldots,v_{2m-2}v_{2m-1}\},$$
where $S$ is any graph on vertex set $\{v_{2a-1},v_{2a},v_{2a+1}\}$ with one or two edges and $1 \leq 2a+1 \leq n-1$.
\end{itemize}
\end{theorem}

Note that in Theorem~\ref{th:K4-}, if $n$ is even and $S$ is a path of length~$2$, then we can assume that its middle point is $v_{2a-1}$.

\begin{figure}[H]
\centering
\includegraphics[scale=0.8]{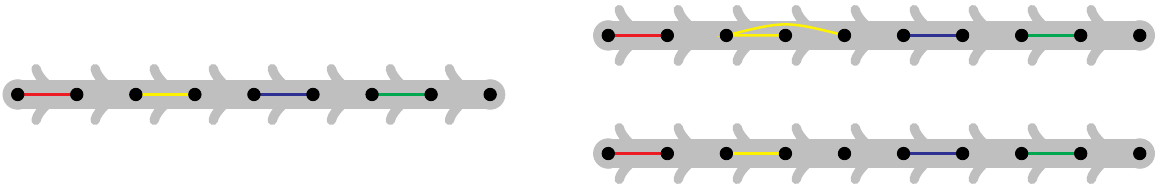}
    \caption{Extremal star-colourings for $K_4^-$. On the left, $n=9$; on the right, $n=10$, $a=2$ and the edges of $S$ are yellow.}
\end{figure}

\begin{theorem}\label{th:K5-}
    As $n \rightarrow \infty$,
    $$
    (1+o(1)) \cdot \left(\tfrac{n}{2}\right)^{3/2} \leq \sar(n,K_5^-) \leq (1+o(1))16\cdot \left(\tfrac{n}{2}\right)^{3/2}.
    $$
\end{theorem}

Our final main result is a vast generalisation of Theorem~\ref{th:K5-}, where we  obtain the order of magnitude of $\sar(n,G)$ for a variety of graphs $G$.
The \emph{join} $G_1+G_2$ of two graphs $G_1$ and $G_2$ is obtained by taking vertex-disjoint copies of $G_1$ and $G_2$ and adding
 all edges between $G_1$ and $G_2$.
The join $H$ of
 two nontrivial forests
is a natural example of a graph with $\va(H)=2$.
Note that $K_5^-$ is the join of the (unique) trees on two and three vertices respectively, and therefore Theorem~\ref{th:K5-} is a special case of a result on $K_2+T$ for a tree $T$ (see Lemma~\ref{lm:K5minus}).
We obtain lower and upper bounds for the star-anti-Ramsey number of the join of two trees each on at least three vertices.
Write $K_{s,t}^-$ for the graph obtained from $K_{s,t}$ by removing one edge.

\begin{theorem}\label{th:treejoin}
Let $T_1$ and $T_2$ be trees on $s,t$ vertices where $3 \leq s \leq t$.
Then there is $c=c(s,t)$ such that
$$
{\rm ex}(n,K_{s,t}^-)/2 \leq \sar(n,T_1 + T_2) \leq c n^{2-1/s}.
$$
\end{theorem}
The order of magnitude of ${\rm ex}(n,K_{s,t}^-)$ is known to be $n^{2-1/s}$ when $t$ is much larger than $s$.
The lower bounds for both Theorems~\ref{th:K5-} and~\ref{th:treejoin} are obtained by taking a lexically coloured complete graph and adding a sparse rainbow graph (as in Theorem~\ref{th:K4-}).

\subsubsection{An open problem}

A major open problem in extremal graph theory is to determine for which $r \in [1,2]$ there is a graph $H$ such that ${\rm ex}(n,H) = \Theta(n^r)$. We say that such an $r$ is \emph{extremal realisable}.

Bukh and Conlon~\cite{BukhConlon}, solving a weaker but also longstanding problem, showed that, if we replace $H$ by a family $\mc{H}$ of graphs, then every rational $r \in [1,2]$ is extremal realisable.
Since then, many $r$ have been shown to be extremal realisable (by a single graph $H$), including 
$1+p/q$ for positive integers $p,q$ with $q>p^2$~\cite{JiangQiu}, and $2-a/b$ for $b \geq \max\{a,(a-1)^2\}$~\cite{ConlonJanzer}.

We can ask analogous questions in the context of star-colourings.     We say that $r$ is \emph{star realisable} if there is a graph $H$ with $\sar(n,H) = \Theta(n^r)$.
 
\begin{problem}
   Which $r \in [1,2]$ are star realisable?
\end{problem}

Theorem~\ref{th:treejoin} combined with existing results about the extremal numbers of very unbalanced complete bipartite graphs yields a countable set of star realisable densities.

\begin{corollary}\label{th:starreal}
    For every integer $s \geq 1$, we have that $2-1/s$ is star realisable.
\end{corollary}

We state several more open problems in Section~\ref{sec:conclude}.

\subsection{Organisation and notation}

We state some results on extremal numbers and some probabilistic tools which will be used in our proofs in Section~\ref{sec:prelim}.
Section~\ref{sec:lower} introduces some special colourings which give (usually sharp) lower bounds.
%We prove the observations Lemmas~\ref{lm:obs1} and~\ref{th:K3} in Section~\ref{sec:lex}.
In Sections~\ref{sec:cycles},~\ref{sec:K4-} and~\ref{sec:K4} we prove Theorems~\ref{th:cycle},~\ref{th:K4-} and~\ref{th:K4} respectively.
Section~\ref{sec:rest} contains the proofs of Theorems~\ref{th:K5-},~\ref{th:treejoin} and Corollary~\ref{th:starreal}. %and~\ref{lm:ea3}. 
We finish with some concluding remarks and open problems in Section~\ref{sec:conclude}.

\medskip
We will use capital letters, e.g.~$G$ to denote a graph, script letters, e.g.~$\ms{G}$ to denote a star-coloured graph, and arrows over capital letters, e.g.~$\ova{G}$ to denote a digraph. Often, $\ova{G}$ will denote an orientation of a graph $G$.

Given a graph $G$, we write $v(G)$ and $e(G)$ for its number of vertices and edges respectively.
Given a subset $A \subseteq V(G)$ and a subset $B \subseteq E(G)$, we write $G-A$ for the graph obtained by removing all vertices in $A$ and all edges incident to any  vertex of $A$,
and $G \sm B$ for the graph with vertex set $V(G)$ and edge set $E(G)\setminus B$.
If $A=\{v\}$ we write $G-v := G-\{v\}$ and if $B=\{e\}$ we write $G \sm e := G\sm\{e\}$.
We write $\overline{A} := V(G) \sm A$.
We make analogous definitions for star-coloured graphs and digraphs.

Given a star-coloured graph $\ms{G}$, we write $C(\ms{G})$ for the set of colours appearing in the colouring. Given an edge $xy$ in $\ms{G}$, we write $c(xy)$ for the colour of $xy$.
We simply write \emph{star} to mean the set of all edges of some colour inside $\ms{G}$. 
We say that a vertex $x$ is a \emph{centre} of a star of colour $c$ if $x$ is a vertex of maximum degree in this star. A \emph{leaf} of a star is a vertex of degree one.
So a single-edge star has two centres and two leaves, while larger stars have a unique centre and every other vertex is a leaf.  
For each vertex $x \in V(\ms{G})$, write $\cen(x)$ for the number of stars centred at $x$, and $\cen_1(x)$ for the number of single-edge stars centred at $x$.

An \emph{orientation induced by $\ms{G}$} is an orientation $\ova{G}$ of the uncoloured graph of $\ms{G}$ where we direct each edge from centre to leaf, where a single centre of a single-edge star is chosen arbitrarily.

\section{Preliminaries}\label{sec:prelim}

\subsection{The Zarankiewicz problem}\label{sec:zaran}

Given integers $s \leq t \leq m, n$, let $z(m,n;s,t)$ be the maximum number of edges in a bipartite graph with parts $A,B$ of sizes $m,n$ respectively such that there are no $S \subseteq A$ and $T \subseteq B$ with $|S|=s$ and $|T|=t$ such that $G[S,T]$ induces a copy of $K_{s,t}$.
Determining $z(m,n;s,t)$ is known as the \emph{Zarankiewicz problem}.
We have
\begin{equation}\label{eq:zex}
2\cdot {\rm ex}(n,K_{s,t}) \leq z(n,n;s,t) \leq 4 \cdot {\rm ex}(n,K_{s,t}).
\end{equation}
Indeed, for the first inequality, given any graph, we can obtain a bipartite graph by taking two copies of its vertex set and adding the pairs of edges $xy$ whenever $x,y$ are in different parts; this preserves $K_{s,t}$-freeness.
For the second, given a $K_{s,t}$-free bipartite graph $G$ with parts $A,B$ each of size $n$, let $a \in A$ and $b \in B$ be such that $ab \notin E(G)$. If $n$ is even, take an equipartition $A=A_1 \cup A_2$ and similarly for $B$. If $n$ is odd, take an equipartition $A\sm\{a\}=A_1' \cup A_2'$ and let $A_i:=A_i' \cup \{a\}$ and similarly define $B_i',B_i$ using $b$ for $i=1,2$. Then $G = G[A_1', B_1] \cup G[A_2, B_1'] \cup G[A_1, B_2'] \cup G[A_2', B_2]$ since $ab \notin E(G)$, and each of these four graphs has exactly $n$ vertices.

We need the following results on the Zarankiewicz problem.

\begin{theorem}\label{th:KST}
    Let $2 \leq s \leq t$ be integers.
    \begin{enumerate}
        \item {\rm (K\H{o}v\'ari--S\'os--Tur\'an,~\cite{KovariSosTuran})} There is a constant $c=c(s,t)$ such that for all $n$, we have $z(n,n;s,t) < cn^{2-1/s}$.
        \item {\rm (Koll\'ar--R\'onyai--Szab\'o,~\cite{KollarRonyaiSzabo})} Suppose $t > (s-1)!$. Then there is a constant $c'=c'(s,t)$ such that for all $n$, we have $z(n,n;s,t) > c'n^{2-1/s}$.
        \item {\rm (F\"uredi,~\cite{Furedi})} ${\rm ex}(n,K_{2,t+1}) = \frac{1}{2}\sqrt{t}n^{3/2} + O(n^{4/3})$ for all $t \geq 1$.
    \end{enumerate}
\end{theorem}

\subsection{Probabilistic tools}

The following version of dependent random choice is an easy adaptation of Lemma~2.1 and Theorem~3.1 in~\cite{FoxSudakov}.

\begin{lemma}[Dependent random choice, \cite{FoxSudakov}]\label{lm:drc}
Let $n,s,b$ be positive integers where $n$ is sufficiently large.
Let $c := 2 \max\{a^{1/s},3b/s\}$ 
and let $\ova{G}$ be an $n$-vertex digraph with $e(\ova{G}) \geq cn^{2-1/s}$. %where $c=\max\{a^{1/s},3(a+b)/s\}$.
Then there is a subset $A$ of $V(\ova{G})$ with $|A|=a$ such that all $s$-subsets of $A$ have at least $b$ common outneighbours.
\end{lemma}
\begin{proof}
    Pick $s$ vertices $v_1, \ldots, v_s$ (with repetition) at random and let $A'$ be the common inneighbourhood of all of them.
    Then by linearity of expectation and convexity, $$\mathbb E(|A'|)=\sum_{v \in V(\ova{G})} \left(\frac{d^{+}(v)}{n}\right)^s \ge n\left(\frac{cn^{1-1/s}}{n}\right)^s=c^s.$$
    The expected number of $s$-subsets of vertices that are in $A'$ and have common outneighbourhood of size less than $b$ is at most
    $${n \choose s} \left(\frac{b}{n}\right)^s< \left(\frac{3n}{s}\right)^s\left(\frac{b}{n}\right)^s = \left(\frac{3b}{s}\right)^s.$$
    By linearity of expectation, there exists a choice of $v_1, \ldots, v_s$ such that the  size of this difference is at least
    $c^s-(3b/s)^s \ge a$. Remove one vertex from each $s$-subset of $A'$ that has fewer than $b$ common outneighbours. The resulting set $A$ satisfies the properties of the lemma and has size at least $a$ as required. 
\end{proof}

\section{Lower bounds}\label{sec:lower}

We next describe some important star-colourings of $K_n$ which provide (in some cases, tight) lower bounds for $\sar(n,H)$.

%------------------------------------------------------------------------------------------------------------------------------------------------------------------------------------------------

\subsection{Special star-colourings}

\noindent
\textbf{The orientable colouring $\ms{G}^{\ova{T}}_n$.}
Let $\ova{T}$ be a tournament with vertices $v_1,\ldots,v_n$.
The $\ova{T}$-colouring $\ms{G}^{\ova{T}}_n$ of $K_n$ is obtained by giving the arc $\overrightarrow{v_iv_j}$ the colour $i$. In other words, the $i$-th star has centre $v_i$ and leaves $N^+_{\ova{T}}(v_i)$. 
We have
\begin{equation}\label{eq:numor}
|C(\ms{G}^{\ova{T}}_n)| = n - |\{v \in V(\ova{T}): d^+_{\ova{T}}(v)=0\}|\in\{n-1,n\}.
\end{equation}

\medskip
\noindent
\textbf{The \sG colouring $\ms{G}^{\lex}_n$.}
Take $\ova{T}$ to be the transitive tournament and let $\ms{G}^{\lex}_n := \ms{G}^{\ova{T}}_n$. Then
\begin{equation}\label{eq:numlex}
|C(\ms{G}^{\lex}_n)| = n - 1.
\end{equation}

\medskip
\noindent
\textbf{The rainbow blow-up colouring $\ms{R}_n(\ms{G}_{n_1},\ldots,\ms{G}_{n_\ell})$.}
Given $\ell \in \mathbb{N}$, integers $n_1,\ldots,n_\ell$ with $n = n_1+\ldots + n_\ell$ and $|n_i-n_j| \leq 1$ for all $i,j \in [\ell]$, and colourings $\ms{G}_{n_i}$ of $K_{n_i}$, take a rainbow colouring with new colours of $T_\ell(n)$, the $n$-vertex Tur\'an graph with $\ell$ parts whose $i$-th part has size $n_i$, and for each $1 \leq i \leq \ell$, add $\ms{G}_{n_i}$.
Then
\begin{equation}\label{eq:numrb}
|C(\ms{R}_n(\ms{G}_{n_1},\ldots,\ms{G}_{n_\ell}))|
= t_{\ell}(n)+\sum_{i \in [\ell]}|C(\ms{G}_{n_i})|.
\end{equation}

\subsection{Modified colourings}

We will need small modifications of the above colourings.
Given any star-colouring $\ms{G}$ of $K_n$,
and a collection $\mc{S}$ of edge-disjoint stars in $K_n$, the \emph{$\mc{S}$-modified colouring} $\ms{G}(\mc{S})$ of $\ms{G}$ is obtained by using a new colour on all of the edges of each star in $\mc{S}$.
We have
\begin{equation}\label{eq:modified}
    |C(\ms{G}(\mc{S}))| = |\mc{S}| + |\{c \in C(\ms{G}): \text{the star of colour }c\text{ is not a subgraph of } \textstyle\bigcup \mc{S}\}|.
\end{equation}
%Every $L$-modified \sG colouring we consider has $e(L)=o(n^2)$.
We will always choose $\mc{S}$ so that $\bigcup\mc{S}$ is a sparse graph; that is, it has $o(n^2)$ edges.
The extremal colouring stated in Theorem~\ref{th:K4-} is an
$\mc{S}$-modified colouring of a lexically coloured~$K_n$. 

We often consider the following special case.
Given a subgraph $L$ of $K_n$, the \emph{$L$-modified colouring} $\ms{G}(L)$ is the colouring $\ms{G}(\mc{S})$ where $\mc{S}$ is the collection of single-edge stars with $\bigcup\mc{S}=L$.
(So $L$ is coloured rainbow.) We have
\begin{align}
    |C(\ms{G}(L))| &= e(L) + |\{c \in C(\ms{G}): \text{the star of colour }c\text{ is not a subgraph of } L\}| \geq e(L) \label{eq:modifiedL}.
\end{align}

\subsection{The lexical colouring}\label{sec:lex}

We first observe some simple facts about star-colourings (recall that~(i) was stated towards the end of  Section~\ref{sec:background}).

\begin{lemma}\label{lm:obs1}~
Let $n\in\mathbb N$ and let  $H$ be a graph. 
\begin{itemize}
    \item[\rm(i)] Every star-colouring of $K_n$ has at least $n-1$ colours.
    \item[\rm(ii)] For every star-colouring of $K_n$ with $n-1$ colours, all centres of colours are distinct.
    \end{itemize}
\end{lemma}

\begin{proof}
Assume one of~(i) and~(ii) does not hold. Then  either there is a star-colouring of $K_n$ that consists of at most $n-2$ coloured stars, or a star-colouring of $K_n$ with $n-1$ colours where not all colour centres are distinct. In both cases, there are distinct vertices $u,v$ which are not centres of any colour. So the edge $uv$ is not coloured, a contradiction. % This proves~(i) and~(ii).
\end{proof}

\begin{lemma}\label{lm:lexical}
The \sG colouring of $K_n$
\begin{enumerate}
	\item contains a rainbow copy of each forest on at most $n$ vertices;
	
	\item does not contain a rainbow cycle  and is the unique star-colouring of $K_n$ with $n-1$ colours that has this property. 
\end{enumerate}
\end{lemma}

\begin{proof}
 It suffices to prove (i) for any $n$-vertex tree $T$. Consider a connected ordering \\ $x_1, x_2, \ldots, x_n$ of $V(T)$ (so $T[\{x_1,x_2,\ldots , x_i\}]$ is connected for each $i$) and embed each $x_i$ into $v_{n-i+1}$. This gives a rainbow copy of $T$ in $K_n$.

For (ii), consider any star-colouring $\ms{G}$ of $K_n$ with $n-1$ colours and let $\ova{G}$ be the induced orientation, which is a tournament. Then by Lemma~\ref{lm:obs1}(ii), $\ova{G}$ has a directed cycle if and only if $\ms{G}$ has a rainbow cycle. Moreover, if $\ova{G}$ does not have a directed cycle, then it is the transitive tournament, and so $\ms{G}$ is lexical. This implies (ii).
\end{proof}

We can already deduce the star-anti-Ramsey number of the triangle, which was first determined by Erd\H{o}s, Simonovits and S\'os~\cite{ErdosSimonovitsSos} in 1975.

\begin{lemma}\label{th:K3}
For  $n \ge 3$, every rainbow $K_3$-free star-colouring of $K_n$ is lexical and $\sar(n, K_3) = n-1$.
\end{lemma}

\begin{proof}
This follows immediately  from Theorem~\ref{th:cycle}, or alternatively  from~\eqref{eq:numlex}, Lemma~\ref{lm:obs1}(i) and Lemma \ref{lm:lexical}(ii).
\end{proof}

Next we prove the second observation from Section~\ref{sec:background} about the existence of $\sar(n,H)$.

\begin{lemma}\label{lm:obs2}~
Let  $H$ be a graph. 
For sufficiently large $n$, $\sar(n,H)$ exists if and only if $\va(H) \geq 2$.
\end{lemma}

\begin{proof}
Note that $\va(H)\ge 2$ if and only if $H$ has a cycle.  If $H$ is a forest then by Theorem \ref{th:canonical}, every star-colouring of a sufficiently large complete graph contains a rainbow copy of $H$ and thus $\sar(n,H)$ does not  exist. If $H$ has a cycle, then $\sar(n,H)$ exists by Lemma \ref{lm:lexical}(ii).
\end{proof}

\subsection{A related Ramsey parameter}\label{sec:nsar}

It is helpful to introduce a further extremal parameter. Let
$$
\nsar(H) := \min\{n \in \mathbb{N}: \text{ every star-colouring of }K_n\text{ contains a rainbow }H\}. 
$$
We can use the canonical Ramsey theorem, Theorem~\ref{th:canonical}, to determine when $\sar(H)$ and $\nsar(H)$ exist.

\begin{lemma}\label{lm:exists}
Let $H$ be a graph. The following are equivalent:
\begin{itemize}
\item[{\rm (1)}] $\nsar(H)$ exists,
\item[{\rm (2)}] $H$ is acyclic,
\item[{\rm (3)}] $\sar(n,H)$ does not exist for all sufficiently large $n$.
\end{itemize}
\end{lemma}

\begin{proof}
    Suppose $\nsar(H)$ exists and is equal to $n_0$. Then clearly $\sar(n,H)$ does not exist for~$n>n_0$. Thus (1) implies (3). 
   Moreover, (3) implies (2) by Lemma~\ref{lm:lexical}(ii).
   Finally, (2) implies (1) since by Theorem~\ref{th:canonical}, any  sufficiently large star-coloured $K_n$ contains a rainbow $K_{v(H)}$   or a lexically coloured  $K_{v(H)}$. By Lemma~\ref{lm:lexical}(i), each of these contains a rainbow~$H$.
Thus $\nsar(H)$ exists.
\end{proof}

Later we will require the following bounds on $\nsar(T)$ for a tree $T$. We denote by $P_t$ the path with $t$ edges.

\begin{lemma}\label{lm:nsar}\hfill
\begin{enumerate}
\item For all integers $t \geq 2$ we have $\nsar(P_{t-1})=t$.
\item If $T$ is a tree with $t$ vertices, then $t \leq \nsar(T) \leq 2(t-1)(t+2)$.
\end{enumerate}
\end{lemma}

\begin{proof}
For~(i), consider any star-colouring of a complete graph and take a tournament which is an induced orientation. This tournament has a directed Hamilton path by R\'edei's theorem~\cite{Redei}, which corresponds to a rainbow path.

The lower bound in~(ii) is clear since $\nsar(T) \geq v(T) = t$. 
The upper bound follows from a result of Jamison, Jiang and Ling~\cite{JamisonJiangLing} who showed that, given any trees $S,T$ with $s,t$ vertices respectively, the minimum $n$ for which every edge-coloured $K_n$ contains either a monochromatic $S$ or a rainbow $T$ satisfies $n \leq (s-2)(t-1)(t+2)$.
Since no star-coloured graph has a monochromatic $P_3$, the bound with $s=4$ applies for $\nsar(T)$. 
\end{proof}

%---------------------------------------------------------------------------------------------------
\subsection{Inductive lower bounds}\label{sec:linsep}

The next simple `linear separation lemma' will be useful.
We first extend the definition of star-anti-Ramsey numbers to families of graphs.
Let $\mc{H}$ be a family of graphs, and define
\begin{align*}
\sar(n,\mc{H}) := \max\{s \in \mathbb{N}: \ &\text{there exists an }s\text{-star-coloured }K_n\text{ without a rainbow copy of }H\\ &\text{for any }H \in \mc{H}\}.
\end{align*}
Given a graph $H$, let $\mc{J}^-_H := \{H-v: v \in V(H)\}$.

\begin{lemma}\label{lm:ls}
	For all graphs $H$ and integers $n \geq 2$, we have
    $$
    \sar(n,H)\geq n-1+\sar(n-1,\mc{J}^-_H).
    $$
\end{lemma}
\begin{proof}
	Given a star-colouring of $K_{n-1}$ which does not contain any rainbow $H-v$ we can add a new vertex $w$ to form a $K_n$. Colouring each edge from $w$ with a new colour adds $n-1$ new colours without creating a rainbow $H$, since the restriction to $K_n-w$ would have to contain a rainbow $H-v$.
\end{proof}

This lemma may be compared to the first inequality in~(\ref{eq:exantiramsey}).
It is useful in the following form. 

\begin{corollary}
	For all graphs $H$ and integers $n \geq 2$, if $\sar(n-1,\mc{J}^-_H)\geq c (n-1)+a$, then $\sar(n,H)\geq (c+1)n+(a-1-c)$.
    \hfill\qed
\end{corollary}

Another lower bound is the following.

\begin{lemma}\label{lm:lowerbd2}
    For all graphs $H$ and integers $n$ where $n \geq v(H)$, we have
    $$
    \sar(n,H) \geq \binom{v(H)-1}{2}+(\delta(H)-1)(n-v(H)+1).
    $$
\end{lemma}

\begin{proof}
    Start with a rainbow $K_{v(H)-1}$. At each stage, take an arbitrary partition of the vertex set into $\delta(H)-1$ parts $A_1,\ldots,A_{\delta(H)-1}$, add a new vertex $x$ and add $\delta(H)-1$ colours which each have centre $x$ and where the leaves of the $i$-th colour are $A_i$. By induction, this colouring contains no rainbow copy of $H$.
\end{proof}

\subsection{Lower bounds for $K_4$}\label{sec:K4lower}

Our proof of Theorem~\ref{th:K4} is rather involved, which is perhaps unsurprising due to multiple ways one can arrive at a lower bound, which turns out to be tight.

\begin{corollary}\label{cor:K4}
	  For all integers $n \geq 4$, we have $\sar(n,K_4)\geq  2n -3$.
\end{corollary}

\begin{proof}
    This follows from Lemmas~\ref{th:K3} and~\ref{lm:ls},
    and also follows from Lemma~\ref{lm:lowerbd2}.
\end{proof}

We give two explicit families of colourings with $2n-3$ colours and without a rainbow $K_4$.

First, let $V_1 \cup V_2$ be a vertex partition of $K_n$, and add lexical colourings with disjoint colour sets into each of $V_1$ and $V_2$. Let $v^* \in V_1$ and for each $v_1 \in V_1\sm\{v^*\}$, use one new colour on every $v_1v_2$ with $v_2 \in V_2$, and finally use different unused colours on $v^*v_2$ for every $v_2 \in V_2$.
The number of colours, counted in the order as described, is $(|V_1|-1)+(|V_2|-1)+(|V_1|-1)+|V_2| = 2n-3$.
Consider any rainbow copy of $K_4$ in this colouring.
There is no rainbow triangle in $V_1$ or in $V_2$, so the copy has exactly two vertices in each $V_i$.
Thus each such vertex in $V_1$ must send at least two different colours to $V_2$. 
But only $v^*$ has this property. Thus the copy of $K_4$ is not rainbow. 

Secondly, let $V_1 \cup V_2 \cup V_3$ be a vertex partition of $K_n$. Add lexical colourings with disjoint colour sets into each of $V_1 \cup V_2 \cup V_3$. For each $i \in [3]$ and $v_i \in V_i$, use one new colour on every $v_iv_{i+1}$ with $v_{i+1} \in V_{i+1}$, where all indices are taken modulo $3$. 
The number of colours is $\sum_{i \in [3]}(|V_i|-1) + \sum_{i \in [3]}|V_i| = 2n-3$.
Since for each $i \in [3]$, each vertex of  $V_i$ sends only one colour to $V_{i+1}$ (where we write $V_4 := V_1$), we know that any rainbow copy of $K_4$ that has one or more vertices in $V_i$ has at most one vertex in $V_{i+1}$. So there must be a $V_i$ that contains at least three vertices of our $K_4$, a contradiction since no $V_i$ has a rainbow triangle.

\subsection{Lower bounds from special colourings}

The colourings described above give lower bounds for our graphs $H$ of interest. In many cases, we will prove matching upper bounds.

\begin{lemma}[Lower bounds]\label{lm:lower}
For all graphs $H$ and integers $n$ where $n \geq v(H)$, we have the following lower bounds for $\sar(n,H)$.
    \begin{enumerate}
    \item \emph{(\cite{AxenovichIverson})} If $\va(H) \geq 3$, then $\sar(n,H) \geq t_{\va(H)-1}(n)+n-\va(H)+1$;
\item for all $k \geq 3$, every colouring as described in Theorem~\ref{th:cycle} 
%has $n+\binom{k-2}{2}-1$ colours and
is rainbow $C_k$-free,     
in particular, $\sar(n,C_k) \geq n + \binom{k-2}{2}-1$;
\item $\sar(n,K_{2k+1}) \geq t_{k}(n) + n -k$ for all integers $ k \geq 2$;
\item $\sar(n,K_{2k}) \geq t_{k-1}(n) + n + \lceil \frac{n}{k-1}\rceil-k-1$ for all integers $k \geq 3$; 
\item $\sar(n,K_4^-) \geq \lfloor 3(n-1)/2\rfloor$.
    \end{enumerate}
\end{lemma}

\begin{proof}

    To prove~(i), assume that  $t=\va(H)\ge 3$. Note that no rainbow copy of $H$ is contained in
 the rainbow blow-up colouring of $K_n$ with $\ell=t-1 \geq 2$ parts where each $\ms{G}_{n_i}$ is a lexical colouring.
The claimed number of colours follows from~\eqref{eq:numlex} and~\eqref{eq:numrb}.

For~(ii), 
note that each
such colouring is an $L$-modified orientable colouring $\ms{G}^{\ova{T}}_n$ where $\ova{T}$ has vertex set $A \cup B$ where $|A|=n-k+1$ and $|B|=k-1$ and $\ova{T} \sm \ova{T}[B]$ is $\ova{C_k}$-free, every arc between $A$ and $B$ is directed towards $B$, and $L$ is the clique on $B$. 
Since $\ova{T}[A]$ has no directed cycle of length $k$ and the vertices of $A$ are each the centre of only one colour, the colouring has no rainbow $C_k$ inside $A$.
Thus every rainbow $C_k$ has at least one vertex in $B$.
 So there is a subpath $P=x_1\ldots x_j$ of $C_k$ contained in $A$ such that
the vertex $u$ preceding $x_1$ and the vertex $w$ after $x_j$ are in $B$.
Then $\ova{x_1u},\ova{x_jw} \in E(\ova{T})$
and since $C_k$ is rainbow, $P$ is a directed path in $\ova{T}$, say from $x_1$ to $x_j$. But then both edges incident to $x_1$ in $C_k$ have colour $x_1$, a contradiction.

Next we prove~(iii) and~(iv).
Note that $\va(K_m)=\lceil m/2\rceil$.
Thus~(iii) follows from (i).
So suppose $m=2k$ is even.
Consider the rainbow blow-up colouring with $\ell := k-1$ parts $V_1,\ldots,V_{k-1}$ where $V_1$ is a largest part, and where $\ms{G}_{n_i}$ is the lexical colouring for $i \in [2,k-1]$ and 
$\ms{G}_{n_1}$ a maximal rainbow $K_4$-free colouring.
The claimed lower bound on the number of colours follows from~\eqref{eq:numrb} and~Corollary~\ref{cor:K4}.
Consider any rainbow clique in this colouring. It has at most two vertices in $V_i$ for each $i \in [2,k-1]$ and at most three in $V_1$. Hence the largest rainbow $K_j$ obeys $j\le 2(k-2)+3=2k-1$.  

Finally, we prove~(v).
Consider the $L$-modified lexical colouring $\ms{G}^{\lex}_n$ of $K_n$ with vertex ordering $v_1,\ldots,v_n$ where $L$ is the matching $v_1v_2,v_3v_4,\ldots,v_{2\lfloor\frac{n-1}{2}\rfloor-1}v_{2\lfloor\frac{n-1}{2}\rfloor}$ of size $\lfloor\frac{n-1}{2}\rfloor$. By~\eqref{eq:numlex} and~\eqref{eq:modified}, the number of colours is $n-1+\lfloor\frac{n-1}{2}\rfloor = \lfloor\frac{3(n-1)}{2}\rfloor$.
Consider any  copy of $K_4^-$, with vertices $v_{i_1},\ldots,v_{i_4}$, where $i_1 < i_2<i_3 < i_4$. Note that of the five edges $v_{i_1}v_{i_2}, v_{i_1}v_{i_3}, v_{i_1}v_{i_4}, v_{i_2}v_{i_3}, v_{i_2}v_{i_4}$ 
%of $K_n$ ---taking this out, there is no $K_n$ here really
only one can be in $L$, and thus they use at most three distinct colours. However, all but at most one of them are present in our copy of $K_4^-$, which therefore is not rainbow.
\end{proof}

Modified lexical colourings provide lower bounds for graphs whose edges cannot be decomposed into two forests.

Given integers $n \geq g \geq 3$, let $\mc{C}_{\leq g} := \{C_3,C_4,\ldots,C_g\}$. We have
\begin{align}
\label{eq:mng}
O_g(1) \cdot n^{1+1/(g-1)} & < {\rm ex}(n,\mc{C}_{\leq g}) < n^{1+1/\lfloor g/2\rfloor} 
& & \text{for $n$ sufficiently large,}\\
\label{eq:exgirth}
{\rm ex}(n,\mc{C}_{\leq 2k+1}) & = (1+o(1))\left(\tfrac{n}{2}\right)^{1+1/k}
& & \text{for $k=2,3,5$}.
\end{align}

The lower bound in~\eqref{eq:mng} follows from the first moment method (see Corollary~2.30 in~\cite{FurediSimonovits}, and the upper bound from the Moore bound, due to Alon, Hoory and Linial~\cite{AlonHooryLinial},
while~\eqref{eq:exgirth} follows from work of Reiman~\cite{Reiman}, Benson~\cite{Benson} and Singleton~\cite{Singleton}.
    
We need to define a further notion of arboricity.
Let $\ea(H)$ denote the \emph{edge arboricity} of $H$, which is the minimum number of parts in an edge partition of $H$ such that $H$ induces a forest in each part.
Burr~\cite{Burr} proved that, for all graphs $H$,
$$
\va(H) \leq \ea(H).
$$
Observe that
\begin{equation}\label{eq:ea}
\ea(H) \geq \left\lceil \frac{e(H)}{v(H)-1}\right\rceil.
\end{equation}
The edge arboricity can be arbitrarily large compared to the vertex arboricity. For example, the join $T_1+T_2$ of two trees $T_1,T_2$ each on $t \geq 2$ vertices has $\va(T_1+T_2) = 2$ but $$
\ea(T_1+T_2) \geq \frac{e(T_1+T_2)}{v(T_1+T_2)-1} = \frac{t^2+2t-2}{2t-1} \geq \frac t2+1.
$$

\begin{lemma}\label{lm:modlexical}
Let $H$ be a graph.
\begin{enumerate}
\item
Let $\mc{J}$ be any family of graphs such that for all forests $F$, we have that $H \sm F$ contains some $J \in \mc{J}$ as a subgraph. Then
$$
\sar(n,H) \geq {\rm ex}(n,\mc{J}).
$$
\item If $\ea(H) \geq 3$ and $g$ is the girth of $H$, then $\sar(n,H) \geq {\rm ex}(n,\mc{C}_{\leq g})$.
\end{enumerate}
\end{lemma}

\begin{proof}
For~(i), 
let $\mc{J}$ be as in the statement and let $L$ be an $n$-vertex graph which is $\mc{J}$-free with ${\rm ex}(n,\mc{J})$ edges. Then for any forest $F$, we have that $L$ is $H \sm F$-free. 
Let $\ms{G}$ be an $L$-modified lexical colouring of $K_n$.
Let $\ms{L}$ be the coloured graph on $L$, which is rainbow.
Suppose there is a rainbow copy $\ms{H}$ of $H$ in $\ms{G}$. Then $\ms{F} := \ms{H} \sm \ms{L}$ is a rainbow subgraph of a lexically coloured graph, so $\ms{F}$ is a rainbow copy of a forest $F$ by Lemma~\ref{lm:lexical}(ii).
By construction $\ms{H} \sm \ms{F}$ is a subgraph of $\ms{L}$.
But $L$ is $H \sm F$-free, a contradiction.
Thus~\eqref{eq:modifiedL} implies that $\sar(n,H) \geq e(L)$, as required.

Next we prove~(ii). 
Given any forest $F$, we have that $H \sm F$ contains a cycle since $\ea(H) \geq 3$.
Thus we can take $\mc{J} := \mc{C}_{\leq g}$ and the assertion follows from~(i).
\end{proof}

\section{Cycles: proof of Theorem~\ref{th:cycle}}
\label{sec:cycles}

First recall that a tournament is \emph{strongly connected} if there is a directed path from $x$ to $y$ for every ordered pair of vertices $(x,y)$. Also recall that a strongly connected $n$-vertex tournament has a directed cycle of each length from $3$ to $n$ (this is a result of Moon from 1966~\cite{Moon}). 

\begin{lemma}\label{lm:allcycles}
	For all $n \geq k \geq 3$, if a star-colouring of $K_n$ has a rainbow $C_k$, then it has a rainbow $C_t$ for all $3\leq t \leq k$.
\end{lemma}
\begin{proof}
 Assume there is a rainbow  $C_k$  and let $\ms{G}$ be the star-coloured subgraph induced by its vertices. Orient the edges of $C_k$ cyclically and orient each of the stars with an edge in $C_k$ consistenly with this edge, either all edges towards the centre or all away from it,  while orienting each of the remaining stars of $\ms{G}$ either away from its centre or towards it. The resulting   orientation of $K_k$   is strongly connected.

 So by Moon's result, the directed graph contains directed cycles of all lengths from 3 to $k$, which translate to rainbow cycles in~$\ms{G}$.
\end{proof}

\begin{lemma}\label{lm:ext}
Let $n \geq 3$ be an integer, let $\ms{G}$ be a star-colouring of $K_n$ and let $v\in V(K_n)$.
	If $v$ is the centre of at least two stars and $\ms{G}-v$ has a rainbow $C_{n-1}$, then $\ms{G}$ has a rainbow $C_n$. 
\end{lemma}
\begin{proof}
	Let $\ms{C}$ be a rainbow Hamiltonian cycle in  $\ms{G}-v$ and give the edges in $\ms{C}$ a cyclic orientation $\ova{C}$.  As in the proof of Lemma \ref{lm:allcycles} assign each star with an edge in $\ova{C}$ an orientation consistent with that of $\ova{C}$.  Orient every other star away from the centre, except for one star centred $v$ which is instead oriented towards $v$.
	
	This gives a strong orientation $\ova{G}$ of $\ms{G}$, so $\ova{G}$ has a directed Hamilton cycle, which must correspond to a rainbow cycle in $\ms{G}$ since it cannot have two consecutive edges from one star. 
\end{proof}

\begin{lemma}\label{lm:hamcycle}
Let $n \geq 3$ be an integer and let $\ms{G}$ be a star-coloured $K_n$.
	If every vertex is the centre of at least two stars, then $\ms{G}$ has a rainbow $C_n$. 
\end{lemma}

\begin{proof}
    We proceed by induction on $n$.
	The statement is true for $n=3$, and we let $n \geq 4$ and assume that it is true for graphs with at most $n-1$ vertices.
    Recall that we write $\cen(x)$ for the number of stars centred at a vertex $x$, and $\cen_1(x)$ for the number of single-edge stars centred at~$x$.

    If there is a vertex $x$ so that $\ms{G}-x$ has at least two stars centred at each vertex, then by induction, $\ms{G}-x$ has a rainbow $C_{n-1}$. By Lemma \ref{lm:ext} we are done.   
    Thus we may assume that for every vertex $v$, there is a vertex $u$ such that $uv$ is a single-edge star and $1 \leq \cen_1(u) \leq \cen(u)=2$.

    Let $A$ be the set of vertices $u$ with $1 \leq \cen_1(u) \leq \cen(u)=2$.
    We say that a triangle is \emph{bad} if its edges are single-edge stars and it has at least two vertices in $A$.
    We claim that there is at most one bad triangle. 
    Suppose, for contradiction, that $v_1u_1u'_1$, $v_2u_2u'_2$ are two bad triangles with $u_1,u_1',u_2,u_2' \in A$. 
    By the definition of~$A$, we have $\cen_1( u ) = \cen( u ) = 2$ for $u \in \{ u_1, u_1',u_2,u_2'\}$.
    Thus, for any $i \in \{ 1,2\}$ and all vertices $w \notin \{ v_i,u_i,u'_i\}$, the edge $w u_i$ is in the star centred at~$w$.
    Hence $\{u_1,u'_1\} = \{u_2,u_2'\}$ and so $ v_1 \ne v_2$ implying that $\cen_1(u_1) \ge |\{v_1,v_2, u_1'\}| = 3$, a contradiction. 
    Let $B$ be vertices of the unique bad triangle if it exists, otherwise $B := \emptyset$.

    Form an auxiliary digraph~$\ova{J}$ with vertex set~$V(\ms{G})$ by adding the arc $\ova{uv}$ whenever $u$ is the centre of exactly one star in $\ms{G}-v$; that is, when $uv$ is a single-edge star and $u \in A$.
    We have $d^-_{\ova{J}}(v) = d^-_{\ova{J}}(v,A) \geq 1$ for all $v \in V(\ms{G})$ and $1 \le d^+_{\ova{J}}(u) \le 2 $ for all $u \in A$. 
    
    Suppose that $d^-_{\ova{J}}(v) \geq 2$ for all $v \notin B$.
    Let $L$ be the set of arcs of $\ova{J}$ whose endpoint is not in $B$. Note that $L$ contains no arc whose startpoint is in $A \cap B$. 
    Thus
    \begin{align*}
        2(n - |B|) \le \sum_{v \notin B} d^-_{\ova{J}}(v) \le 
        |L| \le \sum_{u \in A \setminus B} d^+_{\ova{J}}(u) \le 2 |A \setminus B|.
    \end{align*}
    Hence, we have $V(\ms{G}) \setminus B = A \setminus B$ and, for all $v \in V(\ms{G}) \setminus B$, we have $d^+_{\ova{J}}(v) = d^-_{\ova{J}}(v) = 2$. 
    Note that if $\ova{u w} \in E( \ova{J} )$ with $u,w \in A$, then $\ova{w u} \in E( \ova{J})$.
    Together with $B$, we deduce that $|A| \ge n- |B|$,  that $\ms{G}$ contains a $2$-factor~$F$, whose edges are all single-edge stars, and that $\cen_1( u ) = \cen( u ) = 2$ for all $u \in A$.
    Since $n \ge 4$ and $|B| \in \{0,3\}$, there exist $u,v \in A$ such that $uv \notin E(F)$. 
    Without loss of generality, $uv$ is in the star centred at $u$, so $\cen( u ) \ge |N_F(u) \cup \{v\}| = 3$ contradicting the fact that $u \in A$.

    Therefore, there is a vertex $v \notin B$ with $d^-_{\ova{J}}(v)=1$.
    Let $u \in A$ be such that $\ova{uv} \in E(\ova{J})$.
    Let $z$ be a leaf of another star centred at~$v$ such that $uz$ is not a single-edge star.
    Indeed, we can find such a $z$ because by the hypothesis of the lemma, $v$ is the centre of at least one other star.
    Since $u \in A$, there is at most one $x \neq v$ such that $ux$ is a single-edge star, so we are done unless $\cen_1(v) = \cen(v) = 2$.
     Then $v \in A$ and $vu,vx,ux$ are single-edge stars, so $B = \{x,v,u\}$, a contradiction to $v \notin B$.

    We will form $\ms{G}'$, a new star-coloured $K_{n-1}$, by deleting $v$ and colouring $uz$ with a new colour. 
    If $x$ is a vertex such that the number of stars in~$\ms{G}'$ for which it is a centre is at most one, then $\ova{xv} \in E(\ova{J})$ and so $x=u$.
    But $u$ is in two centres since $uz$ is not a single-edge star in~$\ms{G}$.
    Thus $\ms{G}'$ has at least two stars centred at each vertex and by induction has a rainbow $C_{n-1}$.
    If this $C_{n-1}$ does not contain~$uz$, the rainbow~$C_{n-1}$ is a subgraph of $\ms{G}-v$ and we can apply Lemma~\ref{lm:ext}.
    If it does contain $uz$ we replace it with $uv$ and $vz$, which are from distinct stars with centre~$v$ and so do not appear in the rainbow $C_{n-1}$. 
    Thus we obtain a rainbow~$C_n$.
\end{proof}

We are now able to prove Theorem~\ref{th:cycle}.

\begin{proof}[Proof of Theorem~\ref{th:cycle}]
    We will show the following:
    \begin{itemize}
        \item[($\ast$)] for all integers $n \geq k-1 \geq 2$, in every rainbow $C_k$-free star-colouring of $K_n$ with $n+\binom{k-2}{2}-1$  $=n-k+1+{k-1\choose 2}$ colours there exists a vertex set~$B$ of size~$k-1$ such that 
        \begin{enumerate}
            \item every $x \notin B$ is the centre of precisely one star $S_x$;
            \item the leaves of $S_x$ include $B$ for each $x\notin B$; and
            \item  $B$ forms a rainbow~$K_{k-1}$.
        \end{enumerate} 
    \end{itemize}
    The claimed structure now follows. Indeed, let $\ova{T}$ be the orientation of $\ms{G}\sm\ms{G}[B]$, which is unique since each of the at least two vertices in $B$ is a leaf of every star. A $\ova{C}_k$ in $\ova{T}$ corresponds to a rainbow $C_k$ in $\ms{G}\sm\ms{G}[B]$, so clearly $\ova{T}$ is $\ova{C}_k$-free.

    For proving ($\ast$), fix $k$. We may assume that $k \ge 4$ by Lemma~\ref{th:K3}.
    We proceed by induction on $n$.
    If $n = k-1$, then ($\ast$) holds as $K_{k-1}$ has to be a rainbow coloured if it uses $\binom{k-1}{2}$ colours.
    Thus, we may assume that $n \ge k$. 
    Let $\ms{G}$ be a star-coloured $K_n$ with $\sar(n,C_k)$ colours that does not contain a rainbow~$C_k$. 
    By Lemma~\ref{lm:lower}(ii), we have $\sar(n,C_k) \ge n + \binom{k-2}2-1$. 
    If every vertex is the centre of at least two stars, then there is a rainbow~$C_n$ by Lemma~\ref{lm:ext} and hence a rainbow $C_k$ by Lemma~\ref{lm:allcycles}.
    Thus we may assume there is a vertex $v$ which is the centre of at most one star.
    Note that 
    \begin{align*}
         |C ( \ms{G}-v) | =  \sar(n,C_k) -1  \ge n + \binom{k-2}2-2.
    \end{align*}
    On the other hand, since $\ms{G}-v$ does not contain a rainbow~$C_k$, our induction hypothesis implies that 
     \begin{align*}
         |C ( \ms{G}-v) | \le  \sar(n-1,C_k) =  n -1 + \binom{k-2}2-1.
    \end{align*}
    Therefore, we must have equality and deduce  that $\sar(n,C_k)=n+\binom{k-2}{2}-1$ and $v$ is the centre of precisely one star. 
    Again, by our induction hypothesis, there exists a vertex set~$B_v \subseteq V(\ms{G} - v)$ of size~$k-1$ such that every $x \in V(\ms{G}-v) \setminus B_v$ is the centre of precisely one star, each edge $xy$ with $x \in V(\ms{G}-v) \setminus B_v$ and $y \in B_v$ in the star centred at~$x$ and $B_v$ forms a rainbow $K_{k-1}$.
    
    In particular, setting $B:=B_v$, we obtain that (i) and (iii) hold (where for (i) it is helpful to note that edges  $vw$ with $w\notin B$ must have one of the two colours centred at $v$ and  $w$).
  Let $vx$ be an edge of the star centred at $v$, with $x\in B_v$. If all  edges from $v$ to $B_v$ have the same colour as $vx$, we are done, so assume there is $y\in B_v$ such that the colours of $vy$ and $vx$ are distinct. Construct a rainbow $C_k$ by using $vx, vy$ and an $x$-$y$ path that traverses all vertices of $B_v$ and avoids the only edge at $y$ having the same colour as $vy$, a contradiction.
\end{proof}

%--------------------------------------------

%----------------------------------
\section{$K_4^-$: proof of Theorem~\ref{th:K4-}}
\label{sec:K4-}

In this section we determine~$\sar(n,K_4^-)$ and characterise the family of extremal star-colourings.
We prove the following lemma from which Theorem~\ref{th:K4-} follows easily (as we will see below).

\begin{lemma}\label{lm:K4-2}
For all $n\in \mathbb{N}$, $\sar(n, K_4^-) = \lfloor 3(n-1)/2 \rfloor$.
Moreover, when $n \ge 4$, 
every rainbow $K_4^-$-free star-colouring~$\ms{G}$ of $K_n$ with $\lfloor 3(n-1)/2\rfloor$ colours has the following structure.
\begin{itemize}
\item If $n=2m+1$ is odd, then there exists a vertex partition of $V_1, \dots, V_{m}$ such that 
\begin{enumerate}[label = {\rm(\alph*)}]
	\item $|V_i| = 2$ for all $i \in [m-1]$ and $|V_{m}| = 3$; \label{itm:th:K4-:1}
	\item for all $(x,y) \in V_i \times V_j$ with $1 \le i < j \le m$, the edge $xy$ has colour~$x$; \label{itm:th:K4-:2}
	\item for $i \in [m]$, every edge in $V_i$ has a new unique colour. \label{itm:th:K4-:3}
\end{enumerate}
	\item If $n$ is even, then there exists a vertex set $V_0$ such that 
    \begin{enumerate}[label = {\rm(\alph*$'$)}]
	\item 
    $|V_0| \in \{2,3\}$; 
		\item for all $v \in V_0$ and $u \in V(\ms{G}) \setminus V_0$, the edge $vu$ is coloured~$v$; 
        	\item  there is exactly one colour which only appears between vertices of $V_0$; and 
            	\item $|C(\ms{G} - V_0)| = \sar(n - |V_0|, K_4^-)$. 
 \end{enumerate}
\end{itemize}
\end{lemma}

\begin{proof}[Proof of Theorem~\ref{th:K4-}]
%We proceed by induction on~$n$. 
The theorem is true for $n=2$, and for $n=3$ where the only optimal colouring is the rainbow triangle, which has the required structure.

Next suppose $n \geq 4$ is odd.
Apply Lemma~\ref{lm:K4-2}, and order the vertices as $v_1,\dots,v_n$ such that $V_i = \{v_{2i-1},v_{2i}\}$ for all $i \in [(n-3)/2]$, and $V_{(n-1)/2} = \{v_{n-2},v_{n-1},v_n\}$.

Finally, suppose $n \geq 4$ is even.
Lemma~\ref{lm:K4-2} implies that there is a vertex set $V_0$ of size $2$ or~$3$ with the given structure.
If $|V_0|=2$ then we let $V_0 = \{v_1,v_2\}$ and we are done by induction.
So suppose that $|V_0|=3$.
Let $a := 1$ and let $S$ be graph on~$V_0$ induced by edges of the colour $c$ which only appears between vertices of~$V_0$.
We have $c(xy) \in \{c,x,y\}$ for all distinct $x,y \in V_0$.
Since there are either one or two edges of colour $c$, we can order the vertices of $V_0$ as $v_1,v_2,v_3$ so that $c(v_iv_j) \in \{c,v_i\}$ whenever $1 \leq i < j \leq 3$.
Since $\ms{G}-V_0$ has an odd
number of vertices, $\ms{G}$ has the desired structure.
\end{proof}

It remains to prove Lemma~\ref{lm:K4-2}.

\begin{proof}[Proof of Lemma~\ref{lm:K4-2}]
We proceed by induction on~$n$. 
When $n \le 3$, it is easy to see that $\sar(n, K_4^-) = \lfloor 3(n-1)/2 \rfloor$ and that every colouring has the required form.
Let $n \geq 4$. 
By Lemma~\ref{lm:lower}(v), it suffices to show that $\sar(n, K_4^-) \le \lfloor 3(n-1)/2 \rfloor$ and that every colouring has the required form.
Consider a star-colouring $\ms{G}$ of $K_n$ without a rainbow~$K_4^-$ with $\sar(n, K_4^-)$ colours. 
Note that for any vertex $x$, we have
\begin{align}
\label{eq:centre1}
\sar(n, K_4^-) &= \cen(x)+|C(\ms{G}-x)| 
	\le \cen(x)+\sar(n-1, K_4^-) 
	 = \cen(x)+\lfloor 3(n-2)/2 \rfloor.
\end{align}

\begin{claim} \label{clm:twostars1}
~
\begin{enumerate}[label={\rm (\roman*)}]
	\item We have $1 \le \cen(x) \leq 2$ for all $x \in V(\ms{G})$. 
Moreover, if $n$ is odd, then $\cen(x) = 2$ for all $x \in V(\ms{G})$. \label{itm:twostars:1}
	\item $\sar(n, K_4^-) \le \lfloor 3n/2 \rfloor$. \label{itm:twostars:2}
\end{enumerate}
\end{claim}

\begin{proofclaim}
 Suppose there is a vertex $x$ with $\cen(x) \ge 3$.
 Then let $y_1,y_2,y_3$ be leaves in three distinct stars with centre~$x$. 
Hence $c(xy_1), c(xy_2), c(xy_3)$ are distinct. 
Since each colour class induces a star, none of $c(xy_1)$, $c(xy_2)$ and $c(xy_3)$ appear between $y_1,y_2,y_3$. 
As $\ms{G}$ has no rainbow~$K_4^-$, $y_1y_2y_3$ must form a monochromatic $K_3$, a contradiction.

Suppose there is a vertex with $\cen(x) \leq 1$. 
Then~\eqref{eq:centre1} implies that $\sar(n,K_4^-) \leq \lfloor 3(n-1)/2\rfloor$ 
with equality if and only if $\cen(x)=1$ and $n$ is even.
Recall that Lemma~\ref{lm:lower}(v) states that $\sar(n, K_4^-) \ge \lfloor 3(n-1)/2 \rfloor$.
Hence, \ref{itm:twostars:1} holds. 
Moreover, if we use $\cen(x) \leq 2$ instead in~\eqref{eq:centre1}, then we have $\sar(n, K_4^-) \le \lfloor 3n/2 \rfloor$ and so \ref{itm:twostars:2} holds.
\end{proofclaim}

\medskip
\noindent
\textbf{Case 1: $\cen(z) = 1$ for some $z \in V(\ms{G})$}.

\medskip
\noindent
By Claim~\ref{clm:twostars1}\ref{itm:twostars:1}, $n$ is even. 
%Let $z \in V(\ms{G})$ be such that $\cen(z) = 1$. 
Using~\eqref{eq:centre1}, we have $|C(\ms{G}-z)| = \sar(n-1, K_4^-)$ (see proof of Claim~\ref{clm:twostars1}).
By our induction hypothesis, we can partition $V(\ms{G} - z)$ into parts $V_1, \dots, V_{(n-2)/2}$ satisfying \ref{itm:th:K4-:1}--\ref{itm:th:K4-:3}.
Let $V_1 = \{ x,y\}$ and $c(xy) = c$. 
We also write $z$ for the colour of the star centred at~$z$.

Suppose first that 
$c(xz) = x$ and $c(yz) = y$. Then we are done by setting $V_0 = \{x,y\}$ as $|C(\ms{G} -V_0)| = |C(\ms{G})| - 3 = \sar(n,K_4^-) -3 = \sar(n-2,K_4^-)$. 

Suppose instead that $c(xz) = c$. Then 
the fact that $\cen(z) = 1$ implies that $c(yz) \in \{z,y\}$.
We must have $c(zv)=z$ for all $v \in V(\ms{G}) \sm \{x,y,z\}$ or else $\ms{G}[\{x,y,z,v\}]$ contains a rainbow $K_4^-$ with missing edge $xz$.
In particular, $y$ is not a centre of the star of colour $z$, 
so $\cen(y) = 1$.
Using~\eqref{eq:centre1}, we deduce that $|C(\ms{G}-\{y\})| = \sar(n-1, K_4^-)$.
Moreover, we can partition $V(\ms{G} -\{y\})$ into $V'_1, \dots, V'_{(n-2)/2}$ satisfying \ref{itm:th:K4-:1}--\ref{itm:th:K4-:3}.
Furthermore, we deduce that $V'_i = V_i$ for $i \in [2, (n-2)/2]$.
We are done by setting $V_0 := \{x,y,z\}$ as $|C(\ms{G} -V_0)| = |C(\ms{G})| - 4 = \sar(n,K_4^-) -4 = \sar(n-3,K_4^-)$. 

Without loss of generality, the remaining case is $c(xz) = z$ and $c(yz) \in \{y,z\}$. 
For all $v \in V(\ms{G}) \setminus \{x,y,z\}$, we must have $c(zv) = z$ or else $\ms{G}[\{x,y,z,v\}]$ contains a rainbow $K_4^-$ with missing edge~$yz$. 
We are again done by setting $V_0 := \{x,y,z\}$.

\medskip
\noindent
\textbf{Case 2: $\cen(x) = 2$ for all $x \in V(\ms{G})$.}

\medskip
\noindent
Let $M$ be the subgraph of~$\ms{G}$ induced by edges with a unique colour in~$\ms{G}$ (and no isolated vertices).
As there is no rainbow $K_4^-$ in $\ms{G}$, there is
\begin{itemize}
\item no $3$-edge path in $M$ 
(as such a path could easily be completed to a rainbow $K_4^-$);%(as the remaining edges use at least two colours); ----> I think "remaining edges" is not 100% clear
\item no $3$-edge star in $M$ (by Claim~\ref{clm:twostars1}\ref{itm:twostars:1});
\item no two vertex-disjoint $2$-edge paths in $M$ (if $x_1x_2x_3,y_1y_2y_3$ are such paths, say with  $x_2$ being a centre of the star containing $x_2y_2$, then there is a rainbow $K_4^-$ on $x_1,x_2,x_3,y_2$ with missing edge $y_2x_1$ or $y_2x_3$).
\end{itemize}
Thus $M$ is a matching plus possibly one two-edge path or triangle.
In particular, 
\begin{align}
	e(M) \leq \lfloor(n-3)/2\rfloor + 3 = \lfloor(n+3)/2\rfloor \label{eqn:Mupper}
\end{align}
 with equality if and only if $M$ is a vertex-disjoint union of a matching of size~$\lfloor (n-3)/2 \rfloor$ and a triangle.

Let $V_1, \dots, V_{\ell}$ be the components of~$M$.
Note that $|V_i| = 2$ for all but at most one $i \in [\ell]$.
For each $v \in V(M)$ with $d_M(v) =1$, we also write $v$ for the unique colour of the star with centre~$v$ that is not in~$M$.

%We now study~$\ms{G}[V_1, \dots, V_{\ell}]$. 
We now study $\ms{G}[V_1, \dots, V_{\ell}]$ to show that there is a relabelling $V_1, \dots, V_{\ell}$ of components with
\begin{enumerate}[label = {\rm(\alph*$''$)}]
	\item $|V_i| = 2$ for all $i \in [\ell-1]$ and $|V_{\ell}| \in\{2, 3\}$; \label{itm:th:K4-:1'}
	\item for all $(x,y) \in V_i \times V_j$ with $1 \le i < j \le \ell$, $c(xy) = x$. \label{itm:th:K4-:2'}
\end{enumerate}
To see this, let $\ova{A}$ be an auxiliary oriented graph with vertex set $[\ell]$ where there is an arc $\ova{ij}$ whenever $c(xy) = x$ for all $(x,y) \in V_i \times V_j$. 
Firstly, if $xx'$ and $yy'$ are two vertex-disjoint edges in~$M$, then there is no rainbow path of length~$3$ in $\ms{G}[\{x,x'\},\{y,y'\}]$ or else $\ms{G}[\{x,x',y,y'\}]$ contains at least $5$ colours and hence contains a rainbow $K_4^-$, a contradiction.
This implies that if $c(xy) = x$ for some $(x,y) \in V_i \times V_j$, then $c(x'y') = x'$ for all $(x',y') \in V_i \times V_j$ with $xx',yy' \in E(M)$. 
In particular, exactly one of $\ova{ij}$ or $\ova{ji}$ is an arc of $\ova{A}$ whenever $|V_i|=|V_j|=2$.
Also, if $w \in V_i$ with $d_M(w) =2$, then $c(xw) = x$ for all $x \in V(M) \setminus V_i$. 
Thus $\ova{A}$ is a tournament where if there is $i$ with $|V_i|=3$, the auxiliary vertex $i$ is a sink.
Finally, if there is a oriented triangle $\ova{ijki}$ in $\ova{A}$, then there is a rainbow triangle~$xyzx$ in~$\ms{G}[V_1, \dots, V_{\ell}]$ with $x \in V_i$, $y \in V_j$ and $z \in V_k$ and $|V_i| =2$. Then $\ms{G}[V_i \cup \{y,z\}]$ contains at least $5$ colours and hence contains a rainbow $K_4^-$, a contradiction. 
Thus $\ova{A}$ is a transitive tournament and we relabel components in the transitive order to achieve~\ref{itm:th:K4-:1'} and~\ref{itm:th:K4-:2'}.

By counting colours at each centre, we have 
\begin{align}
		\sar(n, K_4^-)  = |C(\ms{G})| = \sum_y \cen(y) - e(M) = 2n-e(M). \label{eqn:uniqueM}
\end{align}
Together with Claim~\ref{clm:twostars1}\ref{itm:twostars:2}, we have $e(M) \ge \lceil n/2 \rceil \ge 3$. 

Suppose that $|V_1| = 3$. 
%Since $e(M) \ge \lceil n/2 \rceil \ge 3$, 
Then \ref{itm:th:K4-:1'} implies that $M$ is a triangle and $n \in \{4,5,6\}$. 
We may assume that each $v \in V(\ms{G}) \setminus V(M)$ sends edges of colour~$v$ to~$V(M)$, or else $\ms{G}[V(M) \cup \{v\}])$ contains a rainbow~$K_4^-$, a contradiction.
If $|V(\ms{G}) \sm V(M)|=1$, since $\cen(v)=2$, there is an edge between $v$ and $V(M)$ of unique colour in $\ms{G}$, contradicting the definition of~$M$.
Otherwise, $|V(\ms{G}) \setminus V(M)| \in \{2,3\}$, since $\cen(v) = 2$ for all $v \in V(\ms{G}) \setminus V(M)$, the graph $\ms{G} \setminus V(M)$ contains an edge of unique colour in $\ms{G}$, contradicting the definition of~$M$. 

Thus we have that $|V_1|=2$.
We claim that removing $V_1$ from $\ms{G}$ removes exactly three colours.
If not, there exists a vertex $z$ that is the centre of a star with leaf set~$V_1$.
Note that $z \notin V(M)$ by~\ref{itm:th:K4-:2'}.
Let $x \in V_1$ and $w \in V_2$, which exists as $e(M) \ge 3$.
By~\ref{itm:th:K4-:2'}, $\ms{G}[ V_1 \cup \{z,w\}]$ contains a rainbow~$K_4^-$ with missing edge~$xz$, a contradiction. 
This proves the claim.

Therefore, we have 
$\sar(n, K_4^-) = |C(\ms{G})| = 3 + |C(\ms{G}-V_1)| \leq 3+\lfloor 3(n-3)/2\rfloor \leq \lfloor 3(n-1)/2\rfloor$.
Moreover, \eqref{eqn:uniqueM} implies that $|M| \ge \lceil (n+3)/2\rceil$.
By~\eqref{eqn:Mupper}, we have $n$ is odd and $|M| =  (n+3)/2$.
By the remark below~\eqref{eqn:Mupper}, $M$ is a vertex-disjoint union of a matching of size $\lfloor (n-3)/2 \rfloor$ and a triangle.
Together with~\ref{itm:th:K4-:1'} and~\ref{itm:th:K4-:2'},  $V_1, \dots, V_{\ell}$ is the desired partition satisfying \ref{itm:th:K4-:1}--\ref{itm:th:K4-:3}. 
\end{proof}

%%%%%%%%%%%%%%%%%%%%%%%%%%%%%%%%%%%%%%%%%%%%%%%%%%%%%%%%%%%%%%%%%%%%%%%%%%%%%%%%%%%%%%

\section{$K_4$: proof of Theorem~\ref{th:K4}}\label{sec:K4}

Recall that in Section~\ref{sec:K4lower}, we described a large family of star-colourings of $K_n$ with $2n-3$ colours and without a rainbow $K_4$. These multiple extremal colourings suggest that the proof of Theorem~\ref{th:K4} may be quite involved, and indeed it is.

%%%%%%%%%%%%%%%%%%%%%%%%%%%%%%%%%%%%%%%%%%%%%%%%%%%%%%%%%%%%%%%%%%%%%%%%%%%%%%%%%%%%%%%%%
\subsection{Proof ideas}
By Corollary~\ref{cor:K4} (or Lemma~\ref{lm:lower}(iv)), it suffices to show that $\sar(n, K_4) \le 2n-3$.

First we discuss the ideas of the proof. 
We proceed by induction on~$n$. 
Consider a star-colouring~$\ms{G}$ of~$K_n$ without a rainbow~$K_4$ and with $\sar(n, K_4)$ colours. 
For ease of notation, we will write $C(X),C(X,Y)$ respectively for $C(\ms{G}[X]),C(\ms{G}[X,Y])$ whenever $X,Y$ are disjoint vertex sets of~$\ms{G}$.

It is not difficult to show that the statement is true if $\cen(y) \leq 2$ for all vertices $y$, so assume that there is a vertex $x$ with $\cen(x) \geq 3$. 
Let $S_1, \dots, S_k$ be the stars with centre~$x$
and let $V_i := V(S_i)\setminus \{x\}$.
Suppose that all colours 
%in $\bigcup_{i \in [k]}V(S_i)$
on edges between distinct $S_i$ 
also 
appear inside $\ms{G}[V(S_i)]$ for some $i \in [k]$, that is, suppose that
\begin{align}
	C\Big(\bigcup_{i \in [k]}V(S_i)\Big) = \bigcup_{i \in [k]} C(V(S_i)). \label{eqn:ideal1}
\end{align}
Suppose further that 
\begin{align}
	\bigcup_{i \in [k]} V(S_i) = V(K_n). %V(\ms{G}).
    \label{eqn:ideal2}
\end{align}
So 
$ V(K_n) %V(\ms{G}) 
= \{x\} \cup \bigcup_{i \in [k]}V_i$ is a partition.
Using \eqref{eqn:ideal1} and \eqref{eqn:ideal2}, as well as our induction hypothesis which we may use since by assumption $\cen(x) >1$, we have
\begin{align}
	\sar(n,K_4)  &= | C (\ms{G} )| 
	\le \sum_{i \in [k]} |C(V(S_i))|
	\le \sum_{i \in [k]} \sar( v(S_i), K_4)  
 \le \sum_{i \in [k]}(2v(S_i)-3)
	\nonumber
	\\ 
	&= \sum_{i \in [k]}(2|V_i|-1)
	= 2(n-1) -k 
    < 2n-3
	. \label{eqn:ideal3}
\end{align}

In order to prove \eqref{eqn:ideal1}, we need to study the colours appearing between~$S_i$ and~$S_j$ for distinct $i,j \in [k]$.
%We have $|C(x,V_i)| = 1$ for all $i \in [k]$.
Let $\ms{G'} = \ms{G}[V_1, \dots, V_k]$ be the induced $k$-partite subgraph of $\ms{G}$ with vertex classes~$V_1, \dots, V_k$. 
Since $\ms{G}$ is star-coloured and $V_i$ is the set of leaves of star~$S_i$ at~$x$, a rainbow~$K_3$ in~$\ms{G'}$ together with $x$ will form a rainbow~$K_4$.
Thus $\ms{G'}$ is rainbow $K_3$-free. 
One can then show that there is a colour $c_Z$ and vertex subsets $Y, Z \subseteq V(\ms{G})$ such that 
\begin{align*}
	Y \cup Z = \bigcup_{i \in [k]} V(S_i),
	\quad
    | Y \cap Z| = 1
    \quad\text{and}\quad
	C(Y \cup Z) = C(Y) \cup C(Z) \cup \{c_Z\},
\end{align*}
which is the corresponding statement of~\eqref{eqn:ideal1}, see Lemma~\ref{lemma:initialgood}.

Now suppose that \eqref{eqn:ideal2} is false.
Let $W := V(\ms{G}) \setminus \bigcup_{i \in [k]} V_i$, so $x \in W$ and $|W| \ge 2$.
Note that $|W| + |Y|+|Z| = n+2$. 
If $C(\ms{G}) = C(Y \cup Z) \cup C(W)$, then a similar calculation as in~\eqref{eqn:ideal3} also holds.
If $C(\ms{G}) \ne C(Y \cup Z) \cup C(W)$, then we move vertices from~$W$ to~$Y \cup Z$ until they become equal, see Lemmas~\ref{cl:Case2} and~\ref{cl:Case1}.

%%%%%%%%%%%%%%%%%%%%%%%%%%%%%%%%%%%%%%%%%%%%%%%%%%%%%%%%%%%%%%%%%%%%%%%%%%%%%%%%%%%%%%%%%
\subsection{Proof of Theorem~\ref{th:K4}}

We start by formalising the desired cover as highlighted in the previous subsection. 
Let $\ms{G}$ be a star-coloured~$K_n$. 
We say that $\mc{P} = (W,Y,Z,x,v^*,c_Z)$ is \emph{good} if 

\medskip
\begin{enumerate}[label = {\rm P${\arabic*}$:}, ref = {\rm P${\arabic*}$}]
	\item $W, Y, Z \subseteq V(\ms{G})$ and $x,v^* \in V(\ms{G})$ and $c_Z \in C(\ms{G})$; \label{itm:p:1}
	\item $W \cup Y \cup Z =  V(\ms{G})$ and $|W| \geq 1$ and $|Y|,|Z| \ge 2$ and $|W|+|Y|+|Z|=n+2$;\label{itm:p:2}
  \item $C(Y \cup Z) =  B_{\mc{P}}$, where
  $$
  B_{\mc{P}} := C(Y) \cup C(Z) \cup \{c_Z\};
  $$ \label{itm:p:3}
	\item for all $y \in Y \sm \{x\}$ and $z \in Z \sm \{v^*\}$,  the colours in $\ms{G}[\{x\},\{y,z\} \cup W]$ are distinct.\label{itm:p:4}
\end{enumerate}

\medskip
We say that $\mc{P}$ is \emph{great} if $\mc{P}$ also satisfies 

\medskip
\begin{enumerate}[label = {\rm P${\arabic*}$:}, ref = {\rm P${\arabic*}$}, resume]
 
 \item $x \in (W \cap Y) \sm Z$ and $v^* \in (Y \cap Z) \sm W$ and  every vertex $u \in V(\ms{G}) \setminus \{ x,v^*\}$ is in exactly one of $W,Y,Z$;
\label{itm:p:5}
	\item $v^*u$ is in a star centred at~$u$ for all $u \in (Z \setminus \{v^*\}) \cup \{x\}$; \label{itm:p:6}
	\item for all $z \in Z\setminus \{v^*\}$, $xz$ is in a star centred at~$x$ and that star is of colour~$c_Z$.
\label{itm:p:7}
\end{enumerate}

\medskip
Let 
\begin{align*}
	C_{\mc{P}} := C(W) \cup B_{\mc{P}} = C(W) \cup C(Y) \cup C(Z) \cup \{c_Z\}.
\end{align*}
We further say that $\mc{P}$ is \emph{restricted} if $\mc{P}$ is great and, for all $w \in W\sm\{x\}$, we have $C(\{w\},Y) \subseteq B_{\mc{P}} \cup \{c(wx)\}$.

We now state three key lemmas.

\begin{lemma}\label{lemma:initialgood}
Let $\ms{G}$ be a rainbow $K_4$-free star-coloured~$K_n$. 
Let $x$ be a vertex such that $\cen(x) \geq 3$.
Then there exists a great $(W,Y,Z,x,v^*,c_Z)$ such that $xy$ is in a star centred at $x$ for all $y \in Y$. 
\end{lemma}

\begin{lemma}\label{cl:Case2}
Let $\ms{G}$ be a rainbow $K_4$-free star-coloured~$K_n$. 
Suppose that $(W,Y,Z,x,v^*,c_Z)$ is great and $xy$ is in a star centred at $x$ for all $y \in Y$. 
Then there exists a restricted $(W',Y',Z,x,v^*,c_Z)$.
\end{lemma}

\begin{lemma}\label{cl:Case1}
Let $\ms{G}$ be a rainbow $K_4$-free star-coloured~$K_n$. 
Suppose that $(W,Y,Z,x,v^*,c_Z)$ is restricted.
Then there exists a good $\mc{P} = (W',Y,Z',x,v^*,c_Z)$ such that $C(\ms{G}) = C_{\mc{P}}$.
\end{lemma}

We defer their proofs to later subsections. 
We now prove Theorem~\ref{th:K4} using these lemmas. 

\begin{proof}[Proof of Theorem~\ref{th:K4} using Lemmas~\ref{lemma:initialgood}, \ref{cl:Case2} and \ref{cl:Case1}]
By Corollary~\ref{cor:K4} (or Lemma~\ref{lm:lower}(iv)), it suffices to show that $\sar(n, K_4) \le 2n-3$.

We proceed by induction on~$n$.
If $n =2,3$, then the theorem holds. 
Let $n \geq 4$.
Consider a star-colouring $\ms{G}$ of $K_n$ without a rainbow~$K_4$ with $\sar(n, K_4)$ colours. 
Fix a vertex~$x$. 
If $\cen(x) \leq 2$, then 
\begin{align*}
	\sar(n, K_4) \le 2 + |C(\ms{G}-x)| 
	\le 2+ \sar(n-1, K_4) 
	\le 2n-3,
\end{align*}
where the last inequality is due to our induction hypothesis. 
Hence, we may assume that $\cen(x) \ge 3$.
By Lemmas~\ref{lemma:initialgood}, \ref{cl:Case2} and \ref{cl:Case1}, there exists a good $\mc{P} = (W,Y,Z,x,v^*,c_Z)$ such that $C(\ms{G}) = C_{\mc{P}} = C(W) \cup C(Y) \cup C(Z) \cup \{c_Z\}$.
Recall from \ref{itm:p:2} that $|W| \geq 1$ and $|Y|,|Z| \ge 2$ and $|W|+|Y|+|Z|=n+2$. 
In particular, $\max\{|W|,|Y|,|Z|\}<n$.
Then, using the induction hypothesis on $\ms{G}[W]$, $\ms{G}[Y]$ and $\ms{G}[Z]$, we have
\begin{align*}
	\sar(n, K_4) &= |C(\ms{G})| =  |C_{\mc{P}}|
	\le |C(W)| + |C(Y)| +|C(Z)| +1\\
	&\le \sar(|W|, K_4) + \sar(|Y|, K_4) + \sar(|Z|, K_4)+1 \\
	&\le   2|W|-2 + 2|Y|-3 + 2|Z|-3 +1\\
    &= 2n-3,
\end{align*}
where we use $\max\{2|W|-3, 0\}\le 2|W|-2$ in the final inequality.
\end{proof}

%%%%%%%%%%%%%%%%%%%%%%%%%%%%%%%%%%%%%%%%%%%%%%%%%%%%%%%%%%%%%%%%%%%%%%%%%%%%%%%%%%%%%%%%%

\subsection{Proof of Lemma~\ref{lemma:initialgood}}

In this subsection, we prove Lemma~\ref{lemma:initialgood}.
Since $\ms{G}$ is star-coloured without a rainbow $K_4$, the complete partite graph induced by the $\geq 3$ parts which are the leaves of the stars centred at the given vertex $x$ is rainbow $K_3$-free.
The next lemma analyses the structure of such a complete partite graph and shows that it can be almost partitioned into two sets such that all colours appear inside these two sets.
Lemma~\ref{lemma:initialgood} follows easily from this statement.

\begin{lemma} \label{lem:G[V_1V_2]}
Let $\ms{G}$ be a star-colouring of a complete $k$-partite graph $K[V_1,\ldots,V_k]$ with $k \geq 3$. 
Suppose that $\ms{G}$ is rainbow $K_3$-free.
Then there is a partition $Y\cup Z\cup \{v^*\}$ of $V(\ms G)$ such that $Z=V_i$ for some $i$ and
\begin{align*}
	C(\ms{G}[Y, Z]) \subseteq  C( Y \cup \{v^*\} ) \cup C ( Z \cup \{v^*\} ) 
\end{align*}
and every $z \in Z$ is a centre of a star of size at least two which contains $v^*$.
\end{lemma}

First, we show that this lemma implies Lemma~\ref{lemma:initialgood}.

\begin{proof}[Proof of Lemma~\ref{lemma:initialgood} using Lemma~\ref{lem:G[V_1V_2]}]
Let $S_1, \dots, S_k$ be the stars with centre~$x$, so $k \ge 3$. 
For all $i \in [k]$, let $V_i := V(S_i) \setminus \{x\}$
and let $V'  := \bigcup_{i \in [k]} V_i$.
Note that $\ms{G}[V_1,\ldots,V_k]$ has no rainbow~$K_3$ or else there is a rainbow~$K_4$ containing~$x$. 
Thus we can apply Lemma~\ref{lem:G[V_1V_2]} to~$\ms{G}[V_1,\ldots,V_k]$ to obtain a partition $Y'\cup Z'\cup \{v^*\}$ of $V'$ such that  $Z'=V_i$ for some $i$ and
\begin{align}
C(Y',Z') \subseteq C(Y' \cup \{v^*\}) \cup C(Z' \cup \{v^*\}).\label{eqn:YZ}
\end{align}
Moreover, $|C(\{x\},Z')| = 1$  and each $z \in Z'$ is a centre of a star of size at least two and of colour~$c(zv^*)$. Let $c_Z := c(xz)$ for any and thus all $ z \in Z'$. 
Set  
\begin{align*}
	W &:= V(\ms{G}) \sm V', &
	Y &:= Y' \cup \{x,v^*\} = (V' \setminus Z') \cup\{x\}, &
	Z &:= Z' \cup \{v^*\}.
\end{align*}
Note that \ref{itm:p:1}, \ref{itm:p:5}, \ref{itm:p:6} and \ref{itm:p:7} hold by our construction, and \ref{itm:p:2} and \ref{itm:p:3} follow from~\ref{itm:p:5} and~\eqref{eqn:YZ}, respectively. 
To see \ref{itm:p:4}, let $y \in Y \setminus \{x\}$ and $z \in Z \setminus \{ v^*\}$. 
Since $W \cap V' = \emptyset$, for every $w \in W \setminus \{x\}$, the star containing $xw$ has centre $w$, so $c(xw)$ is not equal to $c(xy)$ or $c(xz)$ or $c(xw')$ for any $w' \in W \sm\{w,x\}$. 
Note that $y,z$ are leaves of different stars centred at $x$, so $c(xy) \neq c(xz)$.
This proves~\ref{itm:p:4}.
Hence $(W,Y,Z,x,v^*,c_Z)$ is great. 
Note that $xy$ is in the star centred at $x$ for all $y \in Y$ by construction.
\end{proof}

As in Section~\ref{sec:cycles}, it is helpful to use auxiliary oriented graphs and the next lemma guarantees a structure in an orientation of a complete partite graph which will be useful in proving Lemma~\ref{lem:G[V_1V_2]}.

\begin{lemma} \label{lem:partite}
Let $\ova{G}$ be an orientation of a complete partite graph $K[V_1,\ldots,V_k]$ with $\delta^+(\ova{G}) \ge 1$ that does not contain a copy of $\ova{C_3}$. 
Then there exist vertex classes~$i,i' \in [k]$ 
 such that $\delta^+(\ova{G}[V_i \cup V_{i'}]) \ge 1$, 
 and a directed cycle~$\ova{C}$ in $\ova{G}[V_i \cup V_{i'}]$ such that 
for all $v \in V(\ova{G}) \setminus (V_i \cup V_{i'})$ and $u \in V(\ova{C})$, we have $\ova{vu}\in E(\ova{G})$.
\end{lemma}

\begin{proof}
We recall that $V(\ova{G})$ has a decomposition into \emph{strongly connected components} $U$, where $\ova{G}[U]$ is strongly connected but $\ova{G}[U \cup \{x\}]$ is not strongly connected for any $x \notin U$. For any two strongly connected components $U,U'$, all edges between $U$ and $U'$ are directed from~$U$ to $U'$, or all are directed the other way. This gives rise to an auxiliary oriented graph $\ova{D}$ whose vertices are the strongly connected components and where arcs $\ova{UU'}$ record the common direction of arcs in $\ova{G}$ between components $U$ and $U'$, and $\ova{D}$ is acyclic.
Thus we can choose $U$ such that $U$ is a sink in $\ova{D}$, that is, $U \cap \bigcup_{v \in \overline{U} } N^-_{\ova{G}}(v)  = \emptyset$, and $|U|$ is minimal with this property.

\begin{claim} \label{clm:C4}
Every cycle in $\ova{G}$ intersects exactly two vertex classes.
\end{claim}

\begin{proofclaim}
Suppose there is a cycle $\ova{C} = \ova{v_1v_2\ldots v_\ell}$ intersecting at least three vertex classes, where we consider indices modulo $\ell$, and suppose that it is the shortest such cycle.

Suppose there is a path $\ova{v_{j-1}v_jv_{j+1}}$ in $\ova{C}$ with $v_{j-1},v_j,v_{j+1}$ in distinct vertex classes and at least one other vertex on $\ova{C}$ in the same class as $v_j$. Since $\ova{G}$ is $\ova{C_3}$-free, we have $\ova{v_{j-1}v_{j+1}} \in E(\ova{G})$, and thus we can obtain a cycle by replacing $\ova{v_{j-1}v_jv_{j+1}}$ with $\ova{v_{j-1}v_{j+1}}$.
This cycle meets the same classes as $\ova{C}$ contradicting minimality.
Thus for all $i \in [k]$ with $|V(\ova{C}) \cap V_i| \geq 2$ and $v_j \in V(\ova{C}) \cap V_i$, there is $i' \in [k]$ such that $v_{j-1},v_{j+1} \in V_{i'}$.
Suppose there is such an $i$, so $\ova{C}$ contains at least two vertices from $V_i$.
Then $|V(\ova{C}) \cap V_{i'}| \geq 2$ so $v_{j-1},v_{j+2} \in V_i$.
Continuing round $\ova{C}$ in both directions, we see that it only intersects $V_i$ and $V_{i'}$, a contradiction.

Thus $\ova{C}$ contains at most one vertex in each class. So $\ova{G}[V(\ova{C})]$ is a tournament containing a directed Hamilton cycle, and hence contains a $\ova{C_3}$, a contradiction.
\end{proofclaim}

Let $\ova{C} = \ova{u_1 \dots u_{\ell}}$ be a cycle inside $U$ write $A := V(\ova{C})$ for its vertex set.
Note that such a cycle exists since $\ova{G}[U]$ is strongly connected, and by Claim~\ref{clm:C4}, $\ova{C}$ intersects exactly two vertex classes, $V_i$ and $V_{i'}$.
Since there is no~$\ova{C_3}$, we have, for all $v \in\overline{V_i \cup V_{i'}}$, that $A \subseteq N^+_{\ova{G}}(v)$ or~$A \subseteq N^-_{\ova{G}}(v)$.
Let
$$
V^+ :=N^+_{\ova{G}}(A)\sm (V_i \cup V_{i'})
\quad\text{and}\quad
V^- := N^-_{\ova{G}}(A)\sm (V_i \cup V_{i'}).
$$
Next, we will show that $V^+ = \emptyset$.
By Claim~\ref{clm:C4}, any edges between $V^+$ and $V^-$ must be from $V^-$ to~$V^+$. 
Note that $V^+ \subsetneq U$ as every vertex in $V^+$ is an outneighbour of a vertex in $U$.
Thus we can replace $U$ with a strongly connected component in $\ova{G}[V^+]$ contradicting the choice of~$U$. 
Hence $V^+ = \emptyset$. 
In particular, 
$\ova{vu} \in E(\ova{G})$ for all $v \in \overline{V_i \cup V_{i'}}$ and $u \in A$.

It remains to show that $\delta^+(\ova{G}[V_i \cup V_{i'}]) \geq 1$.
Suppose otherwise, i.e.~suppose there is $x \in V_i \cup V_{i'}$ such that $d^+_{\ova{G}}(x, V_i \cup V_{i'}) = 0$. 
Thus there is a vertex 
$u \in A \cap N^-_{\ova{G}}(x)$.
Since $\delta^+(\ova{G}) \ge 1$, there exists $v \in \overline{V_i \cup V_{i'}}$ with $\ova{xv} \in E(\ova{G})$. 
Then $\ova{uxv}$ is a~$\ova{C_3}$, a contradiction. 
\end{proof}

We now prove Lemma~\ref{lem:G[V_1V_2]}.
To do so, we construct and analyse a digraph $\ova{G}$ with the same vertex set as $\ms{G}$ which captures the structure of the star-colouring.
For this, we use Lemma~\ref{th:K3} which states that every transversal consisting of one vertex from each partite set of $\ms{G}$ has a lexical colouring.

\begin{proof}[Proof of Lemma~\ref{lem:G[V_1V_2]}]
Let $\mathcal{U}$ be the set of all transversal vertex sets in $\ms{G}$; that is, all $U \subseteq V(\ms G)$ with $|U \cap V_i| = 1$ for all $i \in [k]$.
So $\ms{G}[U]$ is a star-coloured $K_k$ for each $U \in \mc{U}$.

We define a $k$-partite digraph~$\ova{G}$ with vertex classes $V_1, \dots, V_k$ as follows. 
For each $U \in \mathcal{U}$, since $\ms{G}[U]$ is rainbow~$K_3$-free, by Lemma~\ref{th:K3}, there exists an ordering $u_1,\dots,u_k$ of~$U$ and colours $c_1, \dots, c_{k-1}$ such that $c(u_iu_j) = c_i$ for all $1 \le i < j \le k$.
Let $\ova{E}_U$ be the collection of arcs $\ova{u_iu_j}$ for all $1 \leq i < j \leq k$ with $i \neq k-1$.
Let the arc set of $\ova{G}$ be the union of the $\ova{E}_U$ over all $U \in \mc{U}$, deleting multiarcs.

\begin{claim}\label{cl:Hprops}
$\ova{G}$ has the following properties:
\begin{enumerate}[label = {\rm(\roman*)}]
	\item for all $\ova{uv} \in E(\ova{G})$, there exists $v' \in N^+_{\ova{G}}(u)$ not in the same vertex class as~$v$ such that $c(u v) = c(uv')$; \label{itm:H2a}
	\item $\ova{G}$ is an oriented graph; \label{itm:H3}
    \item there are $i^*, j^*$ such that if $u,v$ fulfil $\ova{vu},\ova{uv} \notin E(\ova{G})$, then $u,v\in V_{i^*}\cup V_{j^*}$
    and $d^+_{\ova{G}}(u) = d^+_{\ova{G}}(v)=0$; 
    \label{itm:H1}
	\item for each $u \in V(\ms{G})$, there exists a colour $c_u$ such that $c(uv) = c_u$ for all $\ova{uv} \in E(\ova{G})$ and $c_u$ appears on no other arcs of $\ova{G}$; \label{itm:H5}
	\item if $d^+_{\ova{G}}(u) \ge 1$, then $u$ is the centre of a star of colour~$c_u$ of size at least~$2$;\label{itm:H6}
	\item $\ova{G}$ is $\ova{C_3}$-free. \label{itm:H7}
\end{enumerate}
\end{claim}

\begin{proofclaim}
For~\ref{itm:H2a}, let $\ova{uv} \in E(\ova{G})$.
Let $U \in \mc{U}$ be such that $\ova{uv} \in \ova{E}_U$.
So $u,v \in U$. 
By construction, there exists $v' \in U \setminus \{u,v\}$ such that $c(u v) = c(uv')$, implying~\ref{itm:H2a}.

For~\ref{itm:H3}, if $\ova{G}$ is not oriented, say with $\ova{uv}, \ova{vu} \in E(\ova{G})$, then by~\ref{itm:H2a} we can find vertices~$v',u' \in V$ such that $c(uv') = c(uv) = c(vu')$ contradicting the fact that $\ms{G}$ is star-coloured. 
Thus \ref{itm:H3} holds. 

For~\ref{itm:H1}, suppose there are $v_i,v_j$ with $\ova{v_iv_j},\ova{v_jv_i} \notin E(\ova{G})$ with $v_i \in V_i$ and $v_j \in V_j$ and $i,j$ distinct.
We say that $v_iv_j$ is \emph{missing}.
We need to show that all missing edges lie between $V_i$ and $V_j$.
For every $U \in \mc{U}$ which contains $v_i,v_j$, we have that $\{v_i,v_j\}$ 
plays the role of $\{u_{k-1},u_k\}$, so $\ova{uv_i},\ova{uv_j} \in \ova{E}_U$  
for all $u \in U \sm \{v_i,v_j\}$.
So
\begin{equation}\label{eq:missing}
\overline{V_i \cup V_j} \subseteq N^-_{\ova{G}}(v_i)
\text{ whenever }v_iv_j\text{ is missing with }v_i \in V_i\text{ and }v_j \in V_j.
\end{equation}
If~\ref{itm:H1} does not hold, then there are distinct indices $i_1, i_2, i_3$ and an index $i_4\neq i_3$, as well as  vertices $v_{i_\ell}\in V_{i_\ell}$  for $\ell \in [4]$, so that $v_{i_1}v_{i_2}, v_{i_3}v_{i_4}$ are  missing edges. 
By~\eqref{eq:missing} with $\{i,j\}=\{i_1,i_2\}$, we know that $\ova{v_{i_3}v_{i_1}}, \ova{v_{i_3}v_{i_2}}\in E(\ova G)$. Let $\ell\in[2]$ such that $i_4\neq i_\ell$. Then~\eqref{eq:missing}  with $i=i_3$ and $j=i_4$ implies that $\ova{v_{i_\ell}v_{i_3}}\in E(\ova G)$, a contradiction to~\ref{itm:H3}.
 This completes the proof of~\ref{itm:H1}.

To see~\ref{itm:H5}, we first show that
whenever $\ova{uv_1}, \ova{uv_2} \in E(\ova{G})$, we have $c(uv_1) = c(uv_2)$.
Let $\ova{uv_1},\ova{uv_2} \in E(\ova{G})$ and suppose that $v_1 \in V_j$.
By~\ref{itm:H2a}, there exists $v_1' \in N^+_{\ova{G}}(u) \setminus V_j$ with $c(uv'_1) = c(uv_1)$.
If $v'_1 = v_2$, then we are done. 
Otherwise, pick $U \in \mathcal{U}$ such that $u,v_2 \in U$ and $\{v_1,v_1'\} \cap U \ne \emptyset$ 
which is possible since $k \geq 3$.
Then the statement follows by \ref{itm:H3} and our construction.
We have shown that there is a colour $c_u$ such that $c(uv)=c_u$ for all $\ova{uv} \in E(\ova{G})$. Note that there are at least two such edges by~\ref{itm:H2a}. Thus
it remains to show that there are no $u,v,w$ such that $\ova{uv},\ova{wu} \in E(\ova{G})$ and $c(uv)=c(wu)$.
If this holds then~\ref{itm:H2a} applied to $\ova{wu}$ contradicts the fact that $\ms{G}$ is star-coloured.

Part~\ref{itm:H6} follows from~\ref{itm:H2a} and~\ref{itm:H5}. 
Finally, for~\ref{itm:H7},
if $\ova{G}$ contains a~$\ova{C_3}$, then by~\ref{itm:H5} this corresponds to a rainbow~$K_3$ in $\ms{G}$, a contradiction. 
Thus \ref{itm:H7} holds. 
\end{proofclaim}

We claim that, to prove the lemma, it suffices to find $i \in [k]$ and $v^* \in V(\ms{G})$ such that
\begin{enumerate}[label = {\rm(\alph*)}]
\item\label{itm:a1} $v^* \in \overline{V_i}$ and $d^+_{\ova{G}}(v^*) = 0$, 
\item\label{itm:a2} for all $z \in V_i$, we have $N^+_{\ova{G}}(z) \cup N^-_{\ova{G}}(z) = \overline{V_i}$, and
\item\label{itm:a3} for all $y \in \overline{V_i}$ with~$d^+_{\ova{G}}(y) \ge 1$, we have $d^+_{\ova{G}}(y,  \overline{V_i} ) \ge 1$. 
\end{enumerate}
Indeed, set $Y := \overline{V_i}\sm\{v^*\}$ and $Z := V_i$.
Fix $z \in Z = V_i$.
We have $d^+_{\ova{G}}(z) \geq 1$ since $\ova{zv^*} \in E(\ova{G})$ by~\ref{itm:a1} and~\ref{itm:a2}.
By Claim~\ref{cl:Hprops}\ref{itm:H5} and~\ref{itm:H6}, every $u \in N^+_{\ova{G}}(z)$ lies in a common star with centre~$z$, and this star has size at least two.
This proves the last assertion of the lemma. 

Now let $y \in Y = \overline{V_i}\sm\{v^*\}$ and let $c := c(zy)$.
If $\ova{zy} \in E(\ova{G})$, then $y,v^*$ are leaves of a common star centred at $z$, so $c = c(zv^*) \in C(\ms{G}[V_i \cup \{v^*\}]) =  C(\ms{G}[Z \cup \{v^*\}])$ by Claim~\ref{cl:Hprops}\ref{itm:H5}.
Otherwise, \ref{itm:a2} implies that $\ova{yz} \in E(\ova{G})$.
Hence,  $d^+_{\ova{G}}(y,\overline{V_i}) \geq 1$, so by~\ref{itm:a3}, there is
$y' \in \overline{V_i}$ in the outneighbourhood of $y$ along with $z$, so Claim~\ref{cl:Hprops}\ref{itm:H5} implies that they lie in a common star with centre~$y$.
Thus $c=c(yy') \in C(\ms{G}[\overline{V_i}]) = C(\ms{G}[Y \cup \{v^*\}])$.
This completes the proof of the claim.

The rest of the proof will only concern $\ova{G}$, and we will use the properties of $\ova{G}$ deduced in Claim~\ref{cl:Hprops} to find $i,v^*$ which satisfy~\ref{itm:a1}--\ref{itm:a3}.
We split into the cases depending on the base graph of~$\ova{G}$ (by Claim~\ref{cl:Hprops}\ref{itm:H3}, these are the only cases).

%%%%%%%%%%%%%%%%%%%%%%%%%%%%%%%%%%%%%%%
\medskip
\noindent
\textbf{Case 1: $\ova{G}$ is an orientation of a proper subgraph of $K[V_1,\ldots,V_k]$.}

\medskip
\noindent
In this case, $\ova{G}$ has a missing edge $v_1v_2$, and by Claim~\ref{cl:Hprops}\ref{itm:H1}, all such missing edges lie between a common pair of parts, say $V_1$ and $V_2$, and $v_1,v_2$ both have outdegree $0$ in $\ova{G}$.
We are done by setting $v^* := v_1$ and~$i := 3$.
Indeed,~\ref{itm:a1} and~\ref{itm:a2} are clear.
For~\ref{itm:a3}, let $y \in \overline{V_3}$ with $d^+_{\ova{G}}(y)\ge 1$. Then by Claim~\ref{cl:Hprops}\ref{itm:H2a}, we have $d^+_{\ova{G}}(y,\overline{V_3})\ge 1$.

%%%%%%%%%%%%%%%%%%%%%%%%%%%%%%%%%%%%%%
\medskip
\noindent
\textbf{Case 2: $\ova{G}$ is an orientation of $K[V_1,\ldots,V_k]$. }

\medskip
\noindent
Suppose first that $\delta^+(\ova{G}) \ge 1$. 
After possibly relabelling vertex classes, by Lemma~\ref{lem:partite} and Claim~\ref{cl:Hprops}\ref{itm:H7}, there exists a directed cycle~$\ova{C}$ in $\ova{G}[V_1 \cup V_2]$ such that for all $v' \in V(\ova{G}) \sm (V_1 \cup V_2)$ and $u \in V(\ova{C})$, we have $\ova{v'u}\in E(\ova{G})$.
Let $\ova{uv}$ be an arc in~$\ova{C}$. 
Then Claim~\ref{cl:Hprops}\ref{itm:H2a} implies that there exists $v' \in N^+_{\ova{G}}(u) \setminus (V_1 \cup V_2)$, a contradiction to Claim~\ref{cl:Hprops}\ref{itm:H3}. 
Therefore, we have $\delta^+(\ova{G}) = 0$. 

Successively define vertices $v_1,\ldots,v_\ell$ as follows.
For each $i \geq 1$, if there is a vertex $v$ distinct from $v_1,\ldots,v_{i-1}$ with $d^+(v,\ova{G}-\{v_1,\ldots,v_{i-1}\}) = 0$, set $v_i := v$. Otherwise, stop and let $\ell := i-1$.
Since $\delta^+(\ova{G}) = 0$, we have $\ell \ge 1$. 

Suppose first that $v_1, \dots, v_\ell$ do not lie in a common vertex class. 
Let $s$ be the smallest number such that $v_1$ and~$v_s$ lies in different vertex classes, say without loss of generality that $v_1,\ldots,v_{s-1} \in V_{1}$ and $v_s \in V_{2}$. 
Set $v^* := v_1$ and $i :=3$.
We claim that~\ref{itm:a1}--\ref{itm:a3} hold.
Part~\ref{itm:a1} is clear and~\ref{itm:a2} follows from the assumption of Case~2.
To see~\ref{itm:a3}, consider $y \in \overline{V_3}$.
If $y \notin V_{1}$, then $\ova{y v_1} \in E (\ova{G})$ and so $d^+_{\ova{G}}(y,  \overline{V_3} ) \ge d^+_{\ova{G}}(y,  V_1 ) \ge 1$.
If $y \in V_{1} \sm \{v_1,\dots, v_{s-1}\}$, then $\ova{y v_s} \in E (\ova{G})$ and so $d^+_{\ova{G}}(y,  \overline{V_3} ) \ge  d^+_{\ova{G}}(y,  V_2 ) \ge 1$.
Finally if $y \in \{v_1,\dots, v_{s-1}\}$, then $d^+_{\ova{G}}(y) = 0$ by our definitions of $v_1,\dots, v_{s-1}$.

Thus, we may assume that $v_1, \dots, v_\ell$ lie in one vertex class, say~$V_1$. 
Hence $d^+_{\ova{G}}(v_i) = 0$ for all $i \in [\ell]$ and
$\ova{v v_1} \in E(\ova{G})$ for all $v \in \overline{V_1}$.
Moreover, $\ova{G}_\ell := \ova{G}-\{v_1,\ldots,v_{\ell}\}$ has at least $k-1 \geq 2$ vertex classes and $\delta^+(\ova{G}_{\ell}) \ge 1$. 

This together with Claim~\ref{cl:Hprops}\ref{itm:H7} means that we can apply Lemma~\ref{lem:partite} to see that there exist vertex classes~$V'_1$ and~$V'_2$ of~$\ova{G}_{\ell}$ such that $\delta^+(\ova{G}_{\ell}[V_1' \cup V_2']) \geq 1$, and a directed cycle~$\ova{C}$ in $\ova{G}_{\ell}[V'_1 \cup V'_2]$ such that for all $v' \in V(\ova{G}_{\ell}) \setminus (V'_1 \cup V'_2)$ and $u \in V(\ova{C})$, we have $\ova{v'u}\in E(\ova{G}_{\ell})$.
Without loss of generality, we may assume that $V'_2 = V_2$.
We are done by setting $v^* := v_1$ and $i := 3$.
Indeed,~\ref{itm:a1} and~\ref{itm:a2} are again clear by construction.
To see~\ref{itm:a3}, consider $y \in \overline{V_3}$.
If $y \notin V_{1}$, then $d^+_{\ova{G}}(y,  \overline{V_3} ) \ge d^+_{\ova{G}}(y,  \{v_1\} )  = 1$.
If $y \in V_{1} \sm \{v_1,\dots, v_{s-1}\}$, then $d^+_{\ova{G}}(y,  \overline{V_3} ) \ge  d^+_{\ova{G}_\ell}(y,  V_2 \cap V(\ova{C}) ) \ge 1$.
Finally if $y \in \{v_1,\dots, v_{s-1}\}$, then $d^+_{\ova{G}}(y) = 0$ by our definitions of $v_1,\dots, v_{s-1}$. 
\end{proof}

%%%%%%%%%%%%%%%%%%%%%%%%%%%%%%%%%%%%%%%%%%%%%%%%%%%%%%%%%%%%%%%%%%%%%%%%%%%%%%%%%%%%%%%%%

\subsection{Proof of Lemma~\ref{cl:Case2}}

To prove Lemma~\ref{cl:Case2}, we suppose we have a tuple $\mc{P}=(W,Y,Z,x,v^*,c_Z)$ which is great but not restricted.
That is, there exist $w \in W \setminus \{x\}$ and $y \in Y$ such that $c(wy) \notin B_{\mc{P}} \cup \{c(wx)\}$.
It will suffice to show that moving such vertices $w \in W \setminus \{x\}$ into $Y$ yields a tuple $\mc{P}$ that is still great. Then, we can successively define such tuples by moving these vertices, a process which must eventually stop. The resulting tuple is both great and restricted since it has no bad $w$ as above.
The non-trivial part is verifying~\ref{itm:p:3}. We argue that it holds by contradiction by studying the colours between successive bad~$w$, the special vertex~$x$, and vertices $y \in Y$ and $z \in Z$ where the colour $c(yz)$ witnesses the failure of greatness.

\begin{proof}[Proof of Lemma~\ref{cl:Case2}]
Initially, set $W_0 := W$ and $Y_0 := Y$.
Suppose that, for some $i \ge 0$, we have found
$W_0 \supsetneq W_1 \supsetneq \dots \supsetneq W_i$ and $Y_0 \subsetneq Y_1 \subsetneq \dots \subsetneq Y_i$.
We set $\mc{P}_h := (W_h, Y_h,Z, x, v^*, c_Z)$ and $B_h := B_{\mc{P}_h}$ for all $0 \leq h \leq i$.
Consider the following property for $h \geq 0$.

\medskip
\begin{enumerate}[label = {\rm P${\arabic*}(h)$:}, ref = {\rm P${\arabic*}$}]
\setcounter{enumi}{7}
\item if $h \geq 1$, then for all $y_{h} \in Y_{h} \sm Y_{h-1}$, there exists $y_{h-1} \in Y_{h-1} \sm Y_{h-2}$ (where $Y_{-1}=\emptyset$) such that  $c(y_{h}y_{h-1}) \notin B_{{h-1}}  \cup \{c(x y_{h})\}$.\label{itm:p:8}
\end{enumerate}
\medskip
We write ``\ref{itm:p:1}$(h)$ holds'' if $\mc{P}_h$ satisfies~\ref{itm:p:1} and similarly write~\ref{itm:p:2}$(h)$--\ref{itm:p:8}$(h)$.
Suppose that~\ref{itm:p:1}$(h)$--\ref{itm:p:8}$(h)$ hold for all $0 \leq h \leq i$
(that is, $\mc{P}_h$ is great and~\ref{itm:p:8}$(h)$ holds).
Note that~\ref{itm:p:1}$(0)$--\ref{itm:p:8}$(0)$ do indeed hold since $\mc{P}_0$ is great and~\ref{itm:p:8}$(0)$ is vacuous. 
If $\mc{P}_i$ is restricted, then we are done by setting $W' := W_i$ and $Y' := Y_i$.
Thus we may assume that $\mc{P}_i$ is not restricted and so $|W_i| \ge 2$. 
We now construct $\mc{P}_{i+1}$ as follows. 
Let 
\begin{align*}
W^{\rm bad}_i &:= \{ w \in W_i\sm\{x\} : \exists \ y \in Y_i \text{ with }c(wy) \notin B_i \cup \{c(wx)\}\},
\end{align*}
so $W^{\rm bad}_i \ne \emptyset$ since $\mc{P}_i$ is not restricted.
Set 
\begin{align*}
	W_{i+1} &:= W_i \sm W^{\rm bad}_i
    \quad\text{and}\quad
	Y_{i+1} := Y_i \cup W^{\rm bad}_i.
\end{align*}
Note that $1 \le |W_{i+1}| < |W_i|$.
So to prove the lemma, it suffices to show that $\mc{P}_{i+1}$ is great, as the algorithm will terminate. 
Note that $B_0 \subseteq B_1 \subseteq \dots \subseteq B_i$.

Since $Z$ is unchanged, \ref{itm:p:6}$(i+1)$ and \ref{itm:p:7}$(i+1)$ hold. 
By our construction, \ref{itm:p:1}($i+1$), \ref{itm:p:2}($i+1$) and \ref{itm:p:5}($i+1$) hold since $x,v^*$ do not move and since $W_{i+1}\cup Y_{i+1}=W_{i}\cup Y_{i}$.
As $Y_i \subseteq Y_0 \cup W_0$ and $W_i \subseteq W_0$, \ref{itm:p:4}($i+1$) follows from~\ref{itm:p:4}($0$).

Next, we see that~\ref{itm:p:8}$(i+1)$ holds.
Recall that $Y_{i+1} \sm Y_i = W^{\rm bad}_i \neq \emptyset$. Let $y_{i+1} \in Y_{i+1}\sm Y_i$ be arbitrary.
Since $y_{i+1} \in W^{\rm bad}_i$, there is $y_i \in Y_i$ with $c(y_{i+1}y_i) \notin B_{i} \cup \{c(xy_{i+1})\}$.
We need to show that $y_i \in Y_i \sm Y_{i-1} = W^{\rm bad}_{i-1}$.
We are done if $i=0$, so suppose $i \geq 1$.
Suppose for a contradiction that $y_i \in Y_{i-1}$.
Now, $y_{i+1} \notin W^{\rm bad}_{i-1}$, so $c(y_{i+1}y_i) \in B_{i-1} \cup \{c(xy_{i+1})\} \subseteq B_{i} \cup \{c(xy_{i+1})\}$, a contradiction.
Thus~\ref{itm:p:8}($i+1$) holds.

It remains to show that \ref{itm:p:3}($i+1$) holds. 
Suppose to the contrary that there exist $y_{i+1} \in Y_{i+1} \sm Y_i$ and $z \in Z$ such that 
\begin{align}
c(y_{i+1} z ) \notin B_{{i+1}}. \label{eqn:c(y_{i+1} z)}
\end{align}
By \ref{itm:p:8}$(i+1)$--$(1)$, we can find $y_i, \dots, y_0$ such that $y_h\in Y_h\setminus Y_{h-1}$ and $y_{h+1},y_{h}$ satisfy \ref{itm:p:8}$(h)$ for all $0 \le h \le i$. 
Note that they are distinct by definition.

\begin{claim} \label{clm:68}
For all $0 \le h \le i $, we have 
\begin{enumerate}[label = {$(\rm \roman*)$}]
\item\label{itm:P5} $c(y_{h+1}y_h) \notin \{c(y_h z)\} \cup C(\{y_h,y_{h+1}\}, \{ x,y_0, \dots, y_{h-1} \})$;

\item\label{itm:1-10-1} $c(y_{h+1}y_{h}) \notin C(\{y_{h}\}, \{y_{h+2}, \dots, y_{i+1}\})$;

\item\label{itm:i+1z} $c(y_{i+1}z) \notin \{ c(y_{i+1}x)\} \cup C(\{y_{i+1},z\}, \{ y_0, \dots, y_{i} \})$;

\item \label{itm:gen} $\ms{G}[\{x\}, \{z,y_0, \dots, y_{i+1} \}]$ is a rainbow star;

\item \label{itm:c(y_iz)} $c(y_iz) = c(y_i x)$.
\end{enumerate}
\end{claim}

\begin{proofclaim}
First we prove~\ref{itm:P5}.
Since $\{ x,z,y_0, \dots, y_h \} \subseteq Y_h \cup Z$ (using~\ref{itm:p:5}$(h)$ and the definition of~$y_i$), \ref{itm:p:3}($h$) implies that $\{c(y_h u) : u \in \{ x,z,y_0, \dots, y_{h-1} \}\} \subseteq C(Y_h \cup Z) \subseteq B_{h}$.
Since $y_{h+1} \in Y_{h+1} \sm Y_{h} = W^{\rm bad}_{h}$
we have $y_{h+1} \notin W^{\rm bad}_{h-1}$, and $y_0, \dots, y_{h-1} \in Y_{h-1}$, so $C(\{y_{h+1}\}, \{ y_0, \dots, y_{h-1} \}) \in B_{h-1} \cup \{c(y_{h+1}x)\} \subseteq B_{h} \cup \{c(y_{h+1}x)\}$.
In summary, we have 
\begin{align*}
\{c(y_h z)\} \cup C(\{y_h,y_{h+1}\}, \{ x,y_0, \dots, y_{h-1} \}) \subseteq B_{{h}} \cup \{c(y_{h+1}x)\}.
\end{align*}
Thus the definition of~$y_h$ implies~\ref{itm:P5}. 

If \ref{itm:1-10-1} is false, then we have $c(y_{h+1}y_{h}) = c(y_{h}y_{k}) $ for some $h+2 \le k \le i+1$.
Note that $c(x y_k) \ne c(y_{h}y_{k})$ as $y_{h}y_{k}$ is in a star centred at $y_{h}$.
By the definition of~$y_{h}$, we have $c(y_{h}y_{k}) = c(y_{h+1}y_{h}) \notin B_{h} \cup \{c(xy_{h+1})\}$.
Thus $c(y_{h}y_{k}) \notin B_{h} \cup \{c(xy_k)\}$ and so $y_k \in W^{\rm bad}_h \subseteq Y_{h+1}$.
But $y_k \in Y_k \sm Y_{k-1} \subseteq Y_k \sm Y_{h+1}$,
 a contradiction. 
 So \ref{itm:1-10-1} holds.

For \ref{itm:i+1z}, note that since $\{ x,y_0, \dots, y_{i+1} \} \subseteq Y_{i+1}$, we have $C(\{y_{i+1}\}, \{x, y_0, \dots, y_{i} \}) \subseteq C(Y_{i+1}) \subseteq B_{{i+1}}$.
On the other hand, by~\ref{itm:p:3}($i$) and~\ref{itm:p:5}$(i)$, $C(\{z\}, \{x, y_0, \dots, y_{i} \}) \subseteq C(Y_i \cup Z) = B_{i} \subseteq B_{{i+1}}$. 
Thus \ref{itm:i+1z} follows from \eqref{eqn:c(y_{i+1} z)}.

To see~\ref{itm:gen}, first note that $z \neq v^*$ since by~\eqref{eqn:c(y_{i+1} z)} we have $c(y_{i+1}z) \neq B_{i+1}$, while $\{v^*,y_{i+1}\} \subseteq Y_{i+1}$ by~\ref{itm:p:5}($i+1$) so $c(y_{i+1}v^*) \in C(Y_{i+1}) \subseteq B_{i+1}$.
Thus, $z \in Z \sm \{v^*\}$, $y_0 \in Y_0 \sm \{x\}$ (by our choice of~$y_0$) and $\{y_1,\ldots,y_{i+1}\} \subseteq W_0$.
So~\ref{itm:p:4}$(0)$ implies~\ref{itm:gen}.

Suppose that \ref{itm:c(y_iz)} is false, so $c(y_ix) \neq c(y_iz)$. 
We will obtain a contradiction by showing that $\ms{J}_{xz} := \ms{G}[\{y_i,y_{i+1},x,z\}]$ is properly coloured and hence rainbow.
By~\ref{itm:P5} with $h=i$ and since $c(y_ix) \neq c(y_iz)$, we deduce that all edges incident to $y_i$ in $\ms{J}_{xz}$ have distinct colours.
By~\ref{itm:P5} with $h=i$ and~\ref{itm:i+1z} with $h=i$, all edges incident to $y_{i+1}$ in $\ms{J}_{xz}$ have distinct colours.
By~\ref{itm:gen}, all edges incident to $x$ in $\ms{J}_{xz}$ have distinct colours.
Since $z \in Z$, by~\ref{itm:p:7}$(0)$, the edge $xz$ is in a star centred at $x$, so $c(xz) \notin \{c(y_{i+1}z),c(y_iz)\}$.
Together with~\ref{itm:i+1z} with $h=i$, all edges incident to $z$ in $\ms{J}_{xz}$ have distinct colours.
Thus $\ms{J}_{xz}$ is properly coloured, a contradiction.
\end{proofclaim}

Since $y_0 \in Y_0$ is a leaf of a star centred at~$x$ by our assumption in the statement of the lemma, we have $c(y_0x) \neq c(y_0z)$.

Thus, to complete the proof of the lemma, it suffices to obtain the contradiction $c(y_0z)=c(y_0x)$.
For this we proceed by induction on $j$ from above, as follows.
First recall that $c(y_iz)=c(y_ix)$ from Claim~\ref{clm:68}\ref{itm:c(y_iz)}.
We suppose there is $j \le i-1$ such that
\begin{equation}\label{eq:assum}
\begin{cases}
    c(y_{j+1}z) =c(y_{j+1}x) &\mbox{if } j = i-1\\
    c(y_{h-1}z) = c(y_{h-1}x) = c(y_{h+1}y_{h-1}) \text{ for all }j+2 \leq h \leq i &\mbox{if }j \leq i-2
\end{cases}
\end{equation}
and will show that
\begin{equation}\label{eq:goal}
c(y_{j}z) = c(y_{j}x) = c(y_{j+2}y_{j}).
\end{equation}
This induction proves that $c(y_jz)=c(y_jx)=c(y_{j+2}y_j)$ for all $0 \leq j \leq i-1$, which in particular implies the desired contradiction $c(y_0z)=c(y_0x)$.

\begin{claim} \label{clm:69}\hfill
\begin{enumerate}[label = {$(\rm \roman*)$}]
\item\label{itm:zxstar} For all $j+1 \leq h \leq i$ and every $u \neq y_h$, we have that $c(zu) \neq c(y_{h}z)$ and $c(xu) \neq c(y_hx)$;

\item\label{itm:2star} if $j \leq i-2$ then for all $j+2 \leq h \leq i$ and every $u \neq y_{h-1}$, we have that $c(y_{h+1}u) \neq c(y_{h+1}y_{h-1})$. 

\end{enumerate}
\end{claim}

\begin{proofclaim}
Since $c(y_{h}z) = c(y_{h}x)$ by~\eqref{eq:assum}, we have that $y_{h}z$ is in a star centred at $y_{h}$, so~\ref{itm:zxstar} holds.
Since $y_{h+1}y_{h-1}$ is in a star centred at $y_{h-1}$ (whose leaves contain $x$ and $z$) 
by~\eqref{eq:assum}, \ref{itm:2star} holds. 
\end{proofclaim}

Let 
\begin{align*}
\ms{J}_z^j &:= \ms{G}[\{y_j,y_{j+1},y_{j+2},z\}],
&
\ms{J}_x^j &:= \ms{G}[\{y_j,y_{j+1},y_{j+2},x\}],
&
\ms{J}_3^j := \ms{G}[\{y_j,y_{j+1},y_{j+2},y_{j+3}\}],
\end{align*}
where we only consider $\ms{J}_3^j$ when $j \leq i-2$ so that $y_{j+3}$ is defined.
Recall that none of~$\ms{J}_3^j$, $\ms{J}_z^j$, $\ms{J}_x^j$ is properly coloured.
We now consider the following deductions for all $j \leq i-1$:

\begin{enumerate}[label = {\rm(f$_{\arabic*}$)}]

\item\label{itm:i+1,xz} 
%For all $j \leq i-1$, 
All edges incident to $y_{j+1}$ in $\ms{J}_x^j$, in $\ms{J}_z^j$ and (if $j \leq i-2$) in $\ms{J}_3^j$ have distinct colours.
For~$\ms{J}^j_x$, this follows from Claim~\ref{clm:68}\ref{itm:P5} with $h = j$ and $h = j+1$; and then by~\eqref{eq:assum} for $j+1$, the statement holds for $\ms{J}^j_z$ and (if $j \leq i-2$) $\ms{J}_3^j$.

\item\label{itm:3,3} %For all $j \leq i-2$ 
If $j \leq i-2$, 
all edges incident to $y_{j+3}$ in $\ms{J}_3^j$ have distinct colours.
This follows from Claim~\ref{clm:68}\ref{itm:P5} with $h = j+2$ and Claim~\ref{clm:69}\ref{itm:2star} with $(h,u)=(j+2,y_{j})$.

\item\label{itm:j+2,3} %For all $j \leq i-2$, 
If $j \leq i-2$,
all edges incident to $y_{j+2}$ in $\ms{J}_3^j$ have distinct colours.
This follows from Claim~\ref{clm:68}\ref{itm:P5} with $h=j+2,j+1$.

\item\label{itm:j,3} %For all $j \leq i-2$, 
If $j \leq i-2$,
among edges incident to $y_j$ in $\ms{J}_3^j$, only $y_{j+3}y_j$ and $y_{j+2}y_j$ could share a colour.
This follows from Claim~\ref{clm:68}\ref{itm:1-10-1} with $h=j$.

\item\label{itm:3020} %For all $j \leq i-2$, we have 
If $j \leq i-2$,
$c(y_{j+3}y_j) = c(y_{j+2}y_j)$.
This follows from the fact that $\ms{J}^j_3$ is not properly coloured together with
\ref{itm:i+1,xz}, 
%~\ref{itm:i+1,3},
\ref{itm:3,3}, \ref{itm:j+2,3} and~\ref{itm:j,3}.

%%%%%%%%%%%%%%%%%%%%%%%%%%%%%%%%%%%%%%%%%%%%%%%%%%%%

\item\label{itm:z,z} 
%For all $j \leq i-2$, 
All edges incident to $z$ in $\ms{J}_z^j$ have distinct colours.
This follows from Claim~\ref{clm:68}\ref{itm:zxstar} with $h=j+1$ and $u=y_j, y_{j+2}$,
and with $h=j$ and $u=y_{j+2}$.

\item\label{itm:z,j+2} %For all $j \leq i-2$, 
All edges incident to $y_{j+2}$ in $\ms{J}^j_z$ have distinct colours.
If $j = i-1$, then this follows from Claim~\ref{clm:68}\ref{itm:P5} with $h=i$ and Claim~\ref{clm:68}\ref{itm:i+1z}.
If $j \le i-2$, then \ref{itm:3020} implies $y_{j+2}y_j$ is in a star with centre $y_j$ (whose leaves contains~$y_{j+2}$ and~$y_{j+3}$)
 and so $c(y_{j+2}y_j) \notin \{ c(y_{j+2}z), c(y_{j+2}y_{j+1})\}$.
By~\ref{itm:p:8}$(j+2)$, we have $c(y_{j+2}y_{j+1}) \neq c(xy_{j+2})$,
and $c(xy_{j+2})=c(zy_{j+2})$ by~\eqref{eq:assum}.

\item\label{itm:z,j}
%For all $j \leq i-2$, 
Among edges incident to $y_j$ in $\ms{J}^j_z$, only $y_jz$ and $y_{j+2}y_j$ could share a colour.
This follows from Claim~\ref{clm:68}\ref{itm:P5} with $h=j$.

\item\label{itm:0z20} %For all $j \leq i-2$, 
$c(y_jz) = c(y_{j+2}y_j)$.
This follows from the fact that $\ms{J}^j_z$ is not properly coloured together with~\ref{itm:i+1,xz},~\ref{itm:z,z},~\ref{itm:z,j+2} and~\ref{itm:z,j}.

%%%%%%%%%%%%%%%%%%%%%%%%%%%%%%%%%%%%%%%%%%%%%%%%%%%%%%%%%%%%%%%%%%%%%%%%%%%%%%%%%

\item\label{itm:x,x} %For all $j \leq i-1$, 
All edges incident to $x$ in $\ms{J}_x^j$ have distinct colours.
This follows from Claim~\ref{clm:68}\ref{itm:gen}.

\item\label{itm:x,j+2} %For all $j \leq i-2$, 
All edges incident to $y_{j+2}$ in $\ms{J}^j_x$ have distinct colours.
This follows from Claim~\ref{clm:68}\ref{itm:P5} with $h=j+1$ and~\ref{itm:0z20} which implies that $y_{j+2}y_j$ is in a star with centre $y_j$ (and leaf $z$) and hence $c(y_{j+2}y_j) \neq c(y_{j+2}x)$.

\item\label{itm:j,x} %For all $j \leq i-1$, 
Among edges incident to $y_j$ in $\ms{J}_x^j$, only $y_{j}x$ and $y_{j}y_{j+2}$ could share a colour.
This follows from Claim~\ref{clm:68}\ref{itm:P5} with $h=j$ and Claim~\ref{clm:68}\ref{itm:1-10-1} with $h=j$.

\item\label{itm:200x}
%For all $j \leq i-2$, 
$c(y_jx) = c(y_{j+2}y_j)$.
This follows from the fact that $\ms{J}^j_x$ is properly coloured together with~\ref{itm:i+1,xz},~\ref{itm:x,x},~\ref{itm:x,j+2} and~\ref{itm:j,x}.
\end{enumerate}

Finally,~\eqref{eq:goal} follows from \ref{itm:0z20} and~\ref{itm:200x}.
This completes the proof of the lemma.
\end{proof}

%%%%%%%%%%%%%%%%%%%%%%%%%%%%%%%%%%%%%%%%%%%%%%%%%%%%%%%%%%%%%%%%%%%%%%%%%%%%%%%%%%%%%%%%%
\subsection{Proof of Lemma~\ref{cl:Case1}}

To prove Lemma~\ref{cl:Case1}, 
we suppose that $\mc{Q} = (W,Y,Z,x,v^*,c_Z)$ is a restricted tuple with $C_{\mc{Q}} \subsetneq C(\ms{G})$.
Since $\mc{Q}$ is restricted, there are no `bad' colours in $C(W,Y)$ and hence there must be bad colours in $C(W,Z)$ which do not appear elsewhere in~$\ms{G}$. 
We move all $w$ incident to bad colours into $Z$, but this could create new bad colours, so we repeat the process,
and show that eventually we obtain sets $W',Z'$ for which no such colours remain. 
(If at some point $x \in W$ becomes incident to bad colours into $Z$, we may also reallocate the special vertices $x$ and $v^*$ too.)

\begin{proof}[Proof of Lemma~\ref{cl:Case1}]
Initially, set $W_0 := W$,  $Z_0 := Z$, $\mc{P}_0 := (W_0, Y,Z_0, x, v^*, c_Z)$ and $B_0 := B_{\mc{P}_0}$ and $C_0 := C_{\mc{P}_0}$.
If $C_0 = C(\ms{G})$, then we are done (note that $\mc{P}_0$ is restricted and hence good).
Otherwise, we define $\mc{P}_j, B_j, C_j$ as follows. 
For each $j \geq 0$, let
\begin{align*}
W_j^{\rm bad} &:= 
 \{ w \in W^*_j\sm\{x\} : C(\{w\},Z_j) \not\subseteq C(Y) \cup C(Z_j) \cup \{c(wx)\}\},\\
W_{j+1} &:= W_j \setminus W^{\rm bad}_j,\\
   Z_{j+1} &:=
 \begin{cases}
 Z_j \cup W^{\rm bad}_j &\text{if $c(wx) \in C(Z_j \cup \{w\})$ for all $w \in W^{\rm bad}_j$},  \\
 ( Z_j \setminus \{v^*\} ) \cup W^{\rm bad}_j \cup \{x\}
 &\text{otherwise.}
 \end{cases}
 \end{align*}
We set $\mc{P}_j := (W_j, Y,Z_j, x, v^*, c_Z)$ and $B_j := B_{\mc{P}_j}$ and $C_j := C_{\mc{P}_j}$.
We will further show that 
\begin{enumerate}[label = {\rm (\alph*)}]
    \item $\mc{P}_{j+1}$ is good; \label{itm:Case1a}
    \item $B_{j} \subseteq B_{j+1}$; \label{itm:Case1b}
    \item $W_{j+1} \subseteq W_j$ and moreover if $C_j \ne C(\ms{G})$, then $1 \le |W_{j+1}| < |W_j|$. \label{itm:Case1c}
\end{enumerate}
By~\ref{itm:Case1c}, there is some $j^* > 0$ for which $C(\ms{G}) = C_{j^*}$.
Then together with~\ref{itm:Case1a}, we are done by setting $W' := W_j$ and $Z' := Z_j$. 
Thus, to prove the lemma, it suffices to assume that $C_j \ne C(\ms{G})$ and show that \ref{itm:Case1a}--\ref{itm:Case1c} hold for~$\mc{P}_{j+1}$ given that they hold for $\mc{P}_{0},\ldots,\mc{P}_j$.

As before, given $j \geq 0$, we will write e.g.~``\ref{itm:p:1}($j$) holds'' whenever~\ref{itm:p:1} holds for $\mc{P}_j$.
So~\ref{itm:p:1}--\ref{itm:p:7}($0$) hold and we are assuming that~\ref{itm:p:1}--\ref{itm:p:4}($j$) hold.
By~\ref{itm:Case1b} for $\mc{P}_0,\ldots,\mc{P}_j$, we have 
\begin{align}\label{eq:inc}
B_0 \subseteq B_1 \subseteq \ldots \subseteq B_j.
\end{align} 

We start by proving~\ref{itm:Case1c}.
Note that $W_{j+1} \subseteq W_j$ by our construction.
Suppose the second assertion is false, that is, $C_j \ne C(\ms{G})$ and $W^{\rm bad}_j = \emptyset$.
Thus $C(\{w\},Z_j) \subseteq B_j \cup \{c(wx)\}$ for all $w \in W_j$.
By~\ref{itm:p:5}($0$), $x \in W$ and hence $C(W_j,Z_j) \subseteq B_j \cup C(W_j) \subseteq C_j$.
Together with~\ref{itm:p:2}($j$), \ref{itm:p:3}($j$) and~\ref{itm:p:5}($0$), this implies that
\begin{align*}
    C(\ms{G}) & = C(Y) \cup C(W_j) \cup C(Z_j) \cup C(Y,W_j) \cup C(W_j,Z_j) \cup C(Y,Z_j)\\
    & \subseteq C_j \cup C(Y,W_j)  \cup ( B_j \cup C(W_j)) \cup B_j = C_j \cup C(W_j,Y).
\end{align*}
Hence, $C(\ms{G})= C_j \cup C(W_j,Y)$.
Since $C(\ms{G}) \neq C_j$, there are $w \in W_j$ and $y \in Y$ such that $c(wy) \notin C_j$.
If $w \neq x$, then since $\mathcal{P}_0$ is restricted and $W_j \subseteq W$, we obtain the contradiction $c(wy) \in B_0 \cup \{c(wx)\} \subseteq B_j \cup \{c(wx)\}$ from~\eqref{eq:inc}.
If $w=x$, then by~\ref{itm:p:5}($0$), $x \in Y$ so we obtain the contradiction $c(xy) \in C(Y) \subseteq C_j$.
Therefore \ref{itm:Case1c} holds.

We now prove \ref{itm:Case1a} assuming \ref{itm:Case1b} holds.
Note that $x, v^* \in Y$ and $x \in W$ by~\ref{itm:p:5}($0$), so \ref{itm:p:1}($j+1$) holds. 
By our construction,  $|Z_{j+1}|\ge |Z_j|$ and $|W_j| + |Z_j| = |W_{j+1}|+|Z_{j+1}|$ implying~\ref{itm:p:2}($j+1$).
To see \ref{itm:p:4}($j+1$), observe that, for any $y \in Y \setminus \{x\}$ and $z \in Z_j \setminus \{v^*\}$, we have $\{y,z\} \cup W_{j+1} \subseteq \{y,z\} \cup W_{j}$ and so we are done by \ref{itm:p:4}($j$).

It remains to show that~\ref{itm:p:3}($j+1$) holds; that is, $C(Y,Z_{j+1}) \subseteq B_{j+1}$.
Firstly, since $Z_{j+1} \subseteq Z_j \cup W^{\rm bad}_j \cup \{x\}$, we have $C(Y,Z_{j+1}) \subseteq C(Y,Z_j) \cup C(Y,W^{\rm bad}_j \cup \{x\}) \subseteq B_j \cup C(Y,W^{\rm bad}_j)$ where the final inclusion follows from~\ref{itm:p:3}($j$) and~\ref{itm:p:5}($0$) which implies $x \in Y$.
Secondly, since $\mc{P}_0$ is restricted, for any $W' \subseteq W$, we have
$$
C(Y,W') \subseteq C(Y,\{x\}) \cup C(Y,W'\sm\{x\}) \subseteq C(Y) \cup (B_0 \cup C(\{x\},W')) = B_0 \cup C(\{x\},W').
$$
Thirdly, we claim that $C(\{x\},W^{\rm bad}_j) \subseteq C(Z_{j+1})$.
Indeed, if $c(wx) \in C(Z_j \cup \{w\})$ for all $w \in W^{\rm bad}_j$, then $C(\{x\},W^{\rm bad}_j) \subseteq C(Z_j \cup W^{\rm bad}_j) = C(Z_{j+1})$;
while otherwise we have $W^{\rm bad}_j \cup \{x\} \subseteq Z_{j+1}$ so the same conclusion holds.
Putting these three assertions together, we have
\begin{align*}
 C(Y,Z_{j+1}) & \subseteq B_j \cup C(Y,W^{\rm bad}_j) 
	\subseteq  B_j\cup  \left( B_0 \cup C(\{x\},W^{\rm bad}_j) \right) \\
	& \subseteq B_j \cup  B_0 \cup C(Z_{j+1})
    \stackrel{\eqref{eq:inc}}{=} B_j \cup C(Z_{j+1}) \subseteq B_j \cup B_{j+1} \stackrel{\ref{itm:Case1b}}{=} B_{j+1},
\end{align*}
so~\ref{itm:p:3}($j+1$) holds.
Thus \ref{itm:Case1a} holds assuming \ref{itm:Case1b} holds. 

It remains to prove~\ref{itm:Case1b}.
Suppose it is false, so we have $Z_j \not \subseteq Z_{j+1}$.
Thus we have $Z_{j+1} = ( Z_j \sm \{v^*\} ) \cup W^{\rm bad}_j \cup \{x\}$ and there exists $\wt{w} \in W^{\rm bad}_j$ with $c(\wt{w}x) \notin C(Z_j \cup \{\wt{w}\})$.
Since $C(Y) \cup C(Z_j) \not\subseteq C(Y) \cup C(Z_{j+1})$, we have $v^* \in Z_j$ and there exists $\wt{z} \in Z_{j}\sm\{v^*,x\}$ such that $c(\wt{z}v^*) \notin  C(Y) \cup C(Z_{j+1})$.
The vertices $x,v^*,\wt{z},\wt{w}$ are distinct.
Note that $\wt{z},\wt{w} \in Z_{j+1}$, so
\begin{align}
c(\wt{z}v^*) \ne c(\wt{z}\wt{w}). \label{eqn:c(wt{z}v^*)}
\end{align}

Suppose for a contradiction that $\wt{z} \in Z_j \setminus Z$. Then $\wt{z} \in W$.
By~\ref{itm:p:5}($0$) and~\ref{itm:p:6}($0$), for any $u \in Z\sm\{v^*\}$, we have $c(\wt{z}v^*) \neq c(v^*u) \in C(Z)$ and $c(v^*u)=c_Z$.
Since $c(\wt{z}v^*) \notin C(Z_{j+1})$ and $x,\wt{z} \in Z_{j+1}$, we have $c(\wt{z}v^*) \neq c(\wt{z}x)$.
Altogether, $c(\wt{z}v^*) \notin  C(Y) \cup C(Z) \cup \{c_Z\} \cup \{c(\wt{z}x)\} = B_0 \cup \{c(\wt{z}x)\}$. 
However, $v^* \in Y$ contradicting the fact that $\mc{P}_0$ is restricted.
Hence we have 
\begin{align}
	\wt{z} \in Z \setminus \{v^*\}. \label{eqn:wt{z}}
\end{align} 

Let $\ms{J} := \ms{G}[\{x,v^*,\wt{z},\wt{w}\}]$. 
We now show that $\ms{J}$ is a rainbow~$K_4$ which will be a contradiction. 

Since $v^* \in Y$ and $\{x,v^*,\wt{z},\wt{w}\} \subseteq \{x,v^*,\wt{z}\} \cup W_j$, \ref{itm:p:4}($j$) implies that the colours incident to $x$ in $\ms{J}$ are distinct.
By~\eqref{eqn:wt{z}} and~\ref{itm:p:6}($0$), 
the edges $\wt{z}v^*$ and $xv^*$ are in stars centred at $\wt{z},x$ respectively, so
the colours incident to $v^*$ in~$\ms{J}$ are distinct.
Since $\wt{z}x$ is in a star centred at~$x$ by~\ref{itm:p:7}($0$), its colour is not shared with any other edge incident to $\wt{z}$.
Together with~\eqref{eqn:c(wt{z}v^*)}, the colours incident to $\wt{z}$ in $\ms{J}$ are distinct.
Since $c(\wt{w}x) \notin C(Z_j \cup \{\wt{w}\})$ and $v^* \in Z_j$, we have that $\wt{w}x$ does not share a colour with $\wt{w}\wt{z}$ or $\wt{w}v^*$.
Finally, if $c(\wt{w}v^*) = c(\wt{w}\wt{z})$, then there is a star of colour~$c(\wt{w}v^*)$ of size at least~$2$ with centre~$\wt{w}$.
So $c(\wt{w}v^*)$ only appears on edges incident to $\wt{w}$, and not on $\wt{w}x$, contradicting $C(\{\wt{w}\},Y) \subseteq B_{0} \cup \{c(\wt{w}x)\}$ (from $\mc{P}_0$ being restricted).
Thus the colours incident to $\wt{w}$ in $\ms{J}$ are distinct.
Therefore, $\ms{J}$ is a rainbow~$K_4$, a contradiction.
Thus \ref{itm:Case1a}--\ref{itm:Case1c} hold for~$\mc{P}_{j+1}$, which proves the lemma.
\end{proof}

%%%%%%%%%%%%%%%%%%%%%%%%%%%%%%%%%%%%%%%%%%%%%%%%%%%%%%%%%%%%%%%%%%%%%%%%%%%%%%%%%%%%%%%%%

\section{Joins of graphs}
\label{sec:rest}

In this section, we prove Theorems~\ref{th:K5-} and~\ref{th:treejoin}, as well as deriving Corollary~\ref{th:starreal}.

Recall that the \emph{join} $G_1+G_2$ of two graphs $G_1$ and $G_2$ is obtained by taking vertex-disjoint copies of $G_1$ and $G_2$ and adding the edge $x_1x_2$ for every $x_i \in V(G_i)$, $i \in [2]$.
The \emph{internal} edges of $G_1 + G_2$ are the ones induced by $G_1$ or by $G_2$.

In this section, we will obtain bounds for $\sar(n,T_1+T_2)$ for trees $T_1,T_2$.
The upper bounds are in terms of the Zarankiewicz number introduced in Section~\ref{sec:zaran} and the parameter $\nsar(T)$ introduced in Section~\ref{sec:nsar}.
Known results on Zarankiewicz numbers (Section~\ref{sec:zaran}) then imply that for trees $T_1,T_2$ where $s := v(T_1)$ and $T_2$ is significantly larger than $T_1$, these bounds match and we have $\sar(n,T_1+T_2) = \Theta(n^{1-1/s})$.

For the upper bounds, we obtain an uncoloured oriented graph by orienting each star from its centre (choosing a centre arbitrarily for single-edge stars) and keeping one arc per star.
If there are enough colours, this graph is sufficiently dense to find a set of vertices such that any small subset has many common outneighbours. We will use dependent random choice for this to prove Theorem~\ref{th:treejoin}, while for Theorem~\ref{th:K5-} we find a complete bipartite oriented graph.
There is a rainbow copy $\ms{T}_1$ of the smaller tree in this subset, which may reduce the number of common outneighbours due to colours used by this copy.
The larger tree will be embedded in the common outneighbourhood of the vertices of $\ms{T}_1$, and none of the colours appearing in this set have been used before by construction.
Recall that $\nsar(T)$ was defined in Section~\ref{sec:nsar}. 

\begin{lemma}\label{lm:K5minus}
    For every tree $T$ and sufficiently large $n$, we have
    $$
    \sar(n,K_2+T) \leq 4(1+o(1))\cdot z(a,b;2,\nsar(T)+1)
    $$
    where $a,b=(1+o(1))\frac{n}{2}$.
\end{lemma}

\begin{proof}
    Let $m := \nsar(T)$. 
    Let $\eps>0$ and let $n$ be sufficiently large.
    Suppose that $\ms{G}$ is a star-coloured $K_n$ with $4(1+\eps)\cdot z((1+\eps)\tfrac{n}{2},(1+\eps)\tfrac{n}{2};2,m+1)$ colours.
    Orient each star from its centre (choosing a centre arbitrarily in single-edge stars) and keep one arc of each colour, chosen arbitrarily. This gives an oriented graph $\ova G'$.
    Obtain an oriented graph~$\ova{G}$ from~$\ova{G}'$ by independently 
    and uniformly at random assigning each vertex of $V(\ms{G})$ to a part $A$ or $B$ and only keeping arcs oriented from $A$ to $B$.
    We have $\mathbb{E}(|A|)=\mathbb{E}(|B|) = n/2$ and an arc $\overrightarrow{xy}$ of~$\ova{G}'$ is kept independently with probability $1/4$ so $\mathbb{E}(e(\ova{G})) = (1+\eps) \cdot z((1+\eps)\tfrac{n}{2},(1+\eps)\tfrac{n}{2};2,m+1)$ and hence a Chernoff bound (see e.g.~Corollary~2.3 in~\cite{janson2011random}) implies that we may assume that $|A|,|B| < (1+\eps)n/2$ and $e(\ova{G}) > z((1+\eps)\tfrac{n}{2},(1+\eps)\tfrac{n}{2};2,m+1)$, so $e(\ova{G}) > z(|A|,|B|;2,m+1)$.
    Thus there is a $K_{2,m+1}$ (oriented from $A$ to $B$), with vertex partition $\{x_1,x_2\}$ and $\{y_1,\ldots,y_{m+1}\}$.
    Let $\ms{J}$ be the subgraph of $\ms{G}$ induced by this $K_{2,m+1}$.
    By construction, $\ms{J}$ is rainbow and no $y_i$ is a centre of a colour in $\ms{J}$ unless this is the colour of a single-edge star, so none of the colours between the $y_i$ appear in $\ms{J}$.
    If the colour of some $x_iy_j$ is the colour of $x_1x_2$, delete $y_j$; note that there is at most one such $y_j$ since $\ova{G}$ induces a rainbow subgraph of $\ms{G}$. Without loss of generality, $y_j=y_{m+1}$.
    By definition of $\nsar(T)$, $\{y_1,\ldots,y_{m}\}$ contains a rainbow copy of $T$. This copy together with $x_1,x_2$ contains the required copy of $K_2+T$.
\end{proof}

Note that $K_5^- \cong K_2 + P_2$.
We can use Lemma~\ref{lm:K5minus} to prove rather tight bounds for $\sar(K_5^-,n)$.

\begin{proof}[Proof of Theorem~\ref{th:K5-}]
    For the lower bound,
    we have $\ea(K_5^-) \geq 3$ by~\eqref{eq:ea} (in fact there is equality).
    %The girth of $K_5^-$ is five.
    By Lemma~\ref{lm:modlexical}(ii) and~\eqref{eq:exgirth}, we have $\sar(n,K_5^-) \geq {\rm ex}(n,\mc{C}_{\leq 5}) = (1+o(1))\left(\frac{n}{2}\right)^{3/2}$.

    For the upper bound, we have $\nsar(P_2)=3$ by Lemma~\ref{lm:nsar}(i) and ${\rm ex}(m,K_{2,4}) = m^{3/2}+O(m^{4/3})$ by Theorem~\ref{th:KST}(iii).
    By~\eqref{eq:zex}, we have
    $$
    z((1+o(1))\tfrac{n}{2},(1+o(1))\tfrac{n}{2};2,4) \leq 4\cdot{\rm ex}((1+o(1))\tfrac{n}{2},K_{2,4}) \leq 4(1+o(1))\cdot \left(\tfrac{n}{2}\right)^{3/2}.
    $$
    The result now follows from Lemma~\ref{lm:K5minus}.
    \end{proof}

Next, we prove Theorem~\ref{th:treejoin} on the star-anti-Ramsey number of the join $T_1+T_2$ of two trees each on at least three vertices.
The proof of Lemma~\ref{lm:K5minus} generalises from $K_2$ to an arbitrary tree $S$ to give an upper bound of the order $n^{1-1/\nsar(S)}$, which is only of the correct order when $\nsar(S) = v(S)$. Thus we need a slightly different argument in general.

The lower bound comes from a modified lexical colouring. If there is a rainbow copy $\ms{H}$ of $T_1+T_2$, then the non-modified part can only contain a forest subgraph, so most of the graph must lie in the modified part.
The key lemma is the following which we will use to quantify this `most'.

    \begin{lemma}\label{lm:redblue}
        Let $T_1,T_2$ be trees on $s,t$ vertices respectively where $s,t \geq 3$.
        Let $F$ be a subgraph of $T_1+T_2$ which is a forest. Then $(T_1+T_2) \sm F$ 
        either contains an odd cycle or contains $K^-_{s,t}$ as a subgraph.
    \end{lemma}

\begin{proof}[Proof of Theorem~\ref{th:treejoin}]
For the upper bound, let
$$
t' := \max\{\nsar(T_1)+1,\nsar(T_2)+s-1\}
\quad\text{and}\quad
c := \max\{\nsar(T_1)^{1/s},3t'/s\}.
$$

Suppose we have a star-colouring of $K_n$ with $cn^{2-1/s}$ colours.
    Orient each star from its centre (where we choose an arbitrary centre for single-edge stars) and keep one arc of each colour, chosen arbitrarily.  
  Apply Lemma~\ref{lm:drc} with $\nsar(T_1),t'$ playing the roles of $a,b$ to find a set $A$ of size $\nsar(T_1)$ such that every $s$-subset of $A$ has at least $t' \geq \nsar(T_2)+s-1$ common outneighbours.
    By the definition of $\nsar(T_1)$, there is a copy~$T'_1$ of $T_1$ in $A$, with vertex set $X$, so $|X|=s$.
    Let $Y$ be the common outneighbourhood of $X$, so $X \cap Y = \emptyset$ and $|Y| \geq \nsar(T_2)+s-1$. Note that $\ms{G}[X,Y]$ is rainbow. Delete any $y \in Y$ such that there is $xx' \in E(T'_1)$ with the same colour as $\overrightarrow{xy}$ or $\overrightarrow{x'y}$.
    We delete at most one such $y$ for each edge $xx'$, so delete at most $e(T'_1)=s-1$ vertices $y$ in total.
    Thus at least $\nsar(T_2)$ vertices remain in $Y$.
    None of the edges between them use any colours we have already used.
    Thus we can find a rainbow $T_2$ among them which completes the rainbow $T_1+T_2$.

We obtain the lower bound using Lemma~\ref{lm:redblue}.
Let $\mc{J}$ be the family consisting of $K^-_{s,t}$ and all odd cycles.
Let $H := T_1+T_2$ and let $F$ be a forest. 
Lemma~\ref{lm:redblue} implies that $H \sm F$ 
contains some $J \in \mc{J}$.
Lemma~\ref{lm:modlexical}(i) then implies that $\sar(n,H) \geq {\rm ex}(n,\mc{J})$.
Since every graph has a bipartite subgraph containing at least half of its edges, we have ${\rm ex}(n,\mc{J}) \geq {\rm ex}(n,K^-_{s,t})/2$, completing the proof.
\end{proof}

It remains to prove Lemma~\ref{lm:redblue}.

\begin{proof}[Proof of Lemma~\ref{lm:redblue}] 
Colour the edges of $F$ blue and the edges of $T:=(T_1+T_2)\setminus F$ red. For contradiction, assume there is no red odd cycle, no red $K_{s,t}^-$ and no blue cycle. 
We say that $T_i$ is {\it mixed} if it contains both red and blue edges. We claim that for $i=1,2$ the following holds:
\begin{equation}\label{notwomixed}
    \text{If $T_i$ is mixed then $T_{3-i}$ is not.}
\end{equation}
Indeed, assume $T_1$ and $T_2$ are both mixed. Then there are vertices $x_i, y_i, z_i\in V(T_i)$ such that $x_iy_i$ is blue and $y_iz_i$ is red. As there are no blue cycles, we know that one of $x_1y_2$, $x_2y_1$ is red, say  $x_1y_2$. As there are no red triangles, $x_1z_2$ is blue, and thus $y_1z_2$ is red, and hence $y_1y_2, z_1z_2$ are blue. As $x_1y_1y_2x_2$ and $y_1y_2x_2$ do not form blue cycles, we see that $x_1x_2, y_1x_2$ are red. Now $x_1y_2z_2y_1x_2$ is a red odd cycle, a contradiction. This proves~\eqref{notwomixed}.
Next, we show that
\begin{equation}\label{FimpliesF}
    \text{If $T_i$ is blue, then $T_{3-i}$ is blue.}
\end{equation}
For this, assume $T_{3-i}$ contains a red edge $xy$. Observe that  each of $x,y$ sends at most one blue edge to $V(T_i)$, as there are no blue cycles. Since $|V(T_i)|\ge 3$, it follows that there is a vertex $z\in V(T_i)$ such that $xz, yz$ are red. So there is a red triangle $xyz$, a contradiction. This proves~\eqref{FimpliesF}.

Observe that if $T_1$ and $T_2$ are both blue, then at most one edge between $T_1$ and $T_2$ is blue, and thus there is a red $K_{s,t}^-$, which we assumed was not the case. So by~\eqref{FimpliesF},   $T_i$ is not entirely blue for both $i=1,2$. Hence by~\eqref{notwomixed}, there is an $i\in\{1,2\}$ such that $T_i$ is red. Say $i=1$, i.e.~we have
\begin{equation}\label{T1inT}
    T_1 \text{ is red}.
\end{equation}
Root $T_1$ at any vertex and let $L_o$ and $L_e$ be the union of its odd and even levels respectively (where the root belongs to the even levels). Assume $|L_e|\ge |L_o|$, the other case in analogous. Then $|L_e|\ge 2$ and $|L_o|\ge 1$.
Let $S_1\subseteq V(T_2)$ be the set of all those vertices $x$ such that $xy$ is blue for all $y \in L_o$, and set $S_2:= V(T_2)\setminus S_1$. Since there are no red odd cycles, all edges from $S_2$ to $L_e$ are blue. 
As there are no blue cycles and $|L_e|\ge 2$, this means that $|S_2|\le 1$ and therefore $|S_1|\ge 2$. So again using the fact that there are no blue cycles, we see that $|L_o|=1$, and there are no blue edges inside of $S_1$.
We have shown that 
\begin{align}
\label{eq:join1}    
&\text{every edge between $x_o$ and $S_1$ is blue};\\
\label{eq:join2}
&\text{if $S_2=\{y\}$, then every edge from $y$ to $V(T_1) \sm \{x_o\} = L_e$
    is blue and $yx_o$ is red.}
\end{align}

 Suppose $S_1$ contains a red edge $xz$.
 Then one of its endpoints, say $x$, sends only blue edges to $V(T_1)$, because otherwise $x$ and $z$ together with a suitable path in $T_1$ span an odd red cycle, which is not possible.
But then $S_2 = \emptyset$ as otherwise its single vertex lies in a blue $4$-cycle with $x$. Thus the blue graph $F$ is a spanning double-star (one star has centre $x_o$ and leaves $V(T_2) = S_1$, and one star has centre $x$ and leaves $V(T_1) \sm \{x_o\} = L_e$.
Hence there is a red $K_{s,t}^-$, a contradiction.

So we can assume there is no red edge inside $S_1$ and thus
$T_2$ is a star centred at the only vertex $y$ of $S_2$.
By (\ref{eq:join1}) and (\ref{eq:join2}), there is a blue spanning subgraph consisting of two stars between the parts, and hence there is a red $K_{s,t}^-$, a contradiction.
\end{proof}

Finally, we can derive Corollary~\ref{th:starreal}.

\begin{proof}[Proof of Corollary~\ref{th:starreal}]
We know that $1$ is star realisable (e.g.~by $K_3$), 
and $3/2$ is star realisable by $K_5^-$ (Theorem~\ref{th:K5-})
so we need to show that $2-1/s$ is star realisable for all $s \geq 3$.
Take integers $3 \leq s \leq t$ where $t>(s-1)!+1$. Then ${\rm ex}(n,K_{s,t}^-) \geq {\rm ex}(n,K_{s,t-1}) \geq c'n^{2-1/s}$ by Theorem~\ref{th:KST}(ii) where $c'=c'(s,t)$.
Let $T_1,T_2$ be trees with $s,t$ vertices respectively.
Theorem~\ref{th:treejoin} implies that $c'n^{2-1/s} \leq {\rm ex}(n,K_{s,t}^-) \le \sar(n,T_1+T_2) \leq cn^{2-1/s}$ for some constant $c=c(s,t)$. 
Thus $2-1/s$ is star realisable, as required.
\end{proof}

\section{Concluding remarks and open problems}\label{sec:conclude}

\subsection{Specific graphs}
We recall that Theorem~\ref{th:AxenovichIverson} of Axenovich and Iverson determines $\sar(n,H)$ asymptotically when $\va(H) \geq 3$.
We recall that $\va(K_k) \geq 3$ when $k \geq 5$, and that $\sar(n,K_3)$ follows from work of Erd\H{o}s, Simonovits and S\'os~\cite{ErdosSimonovitsSos}, while we determined the exact value of $\sar(n,K_4)$ (see Theorem~\ref{th:K4}).
Lemma~\ref{lm:lower} gives a lower bound for cliques on at least $5$ vertices.

\begin{problem}
    What is the exact value of $\sar(n,K_k)$ for $k \geq 5$?
\end{problem}

Another natural graph to consider is the $3$-dimensional cube $Q_3$, which has $\va(Q_3)=2$.
Lemma~\ref{lm:lowerbd2} implies that $\sar(n,Q_3) \geq \binom{7}{3}+2(n-7) = 2n+21$.
In particular, is it true that $\sar(n,Q_3)$ is linear in $n$?

\subsection{Linear star-anti-Ramsey number}
More generally, one can ask the following.

\begin{problem}
    For which graphs $H$ do we have $\sar(n,H)=\Theta(n)$?
\end{problem}

Given an integer $k \geq 3$ and a tree $T$, there are constants $c_k,c_T>0$ such that for all sufficiently large $n$ we have ${\rm ex}(n,\mc{C}_{\leq 2k}) \geq n^{1+c_k/k}$, and ${\rm ex}(n,T) < c_T n$.
Thus we can characterise those $H$ with linear extremal number: ${\rm ex}(n,H)=\Theta(n)$ if and only if $H$ is a forest with at least two edges.

Schiermeyer and Sot\'ak~\cite{SchiermeyerSotak} proved the corresponding characterisation for the anti-Ramsey problem: if a graph $H$ contains at least two cycles, then $\ar(n,H)$ is superlinear, while if $H$ is unicyclic then $\ar(n,H)=\Theta(n)$.

The corresponding problem for star-colouring seems harder.
Recall from Theorem~\ref{th:AxenovichIverson} and Lemma~\ref{lm:obs1} that to solve the problem, we need to determine which $H$ with $\va(H)=2$ satisfy $\sar(n,H) = O(n)$.
Every graph $H$ with linear star-anti-Ramsey number has edge arboricity at most two:

\begin{lemma}\label{lm:ea3}
    Let $H$ be a graph with $\ea(H) \geq 3$. Then $\sar(n,H) = \Omega(n^{1+1/(v(H)-1)})$.
\end{lemma}

\begin{proof}
This follows from Lemma~\ref{lm:modlexical}(ii),~\eqref{eq:mng}, and the fact that a graph $H$ has girth at most $v(H)$.
\end{proof}

Since we need only consider $H$ with $\va(H) = 2$, this leaves open the question of determining the magnitude of $\sar(n,H)$ when $\ea(H)=2$.

\begin{lemma}\label{lm:star}
Let $H$ be a graph with an edge partition into a star and a forest $F$. Then
$$
n-1 \leq \sar(n,H) \leq (\nsar(F)-1)n.
$$
\end{lemma}

\begin{proof}
Since $2=\ea(H) \geq \va(H)$, Lemma~\ref{lm:obs1}(ii) implies the lower bound.
Let $a := \nsar(F)$.
Suppose there is a vertex which is the centre of $\geq a$ stars and take one arbitrary leaf of each. The graph spanned by these vertices uses disjoint colours from the stars and contains a rainbow copy of $F$, giving the required rainbow copy of $H$.
Thus we may assume that every vertex is the centre of at most $a-1$ stars, so $\sar(n,H) \leq (a-1)n$.
\end{proof}

\begin{problem}\label{prob:ea}
    Let $H$ be a graph with $\ea(H)=2$. Is $\sar(n,H) = \Theta(n)$?
\end{problem}
Here it would be interesting to better understand the case of 4-regular graphs. Recall that by the tree-packing theorem of Nash-Williams \cite{NW61} and Tutte \cite{Tu61}, a well-connected, in the sense defined in the two papers,  
4-regular graph $H$ contains two edge-disjoint spanning trees. Hence the graph has $\ea(H)=3$, with a partition into  two spanning trees and a forest on two edges.   So by deleting two edges from $H$ we can get a graph $H'$  with $\ea(H')=2$. Note that the number of deleted edges does not depend on the size of $H$.  In particular it would be interesting to determine $\sar(n,Q_4)$.

\medskip
\noindent
\textbf{Acknowledgements.}
The work that led to this paper was begun at the workshop \emph{Extremal Graphs and Hypergraphs} in July 2024 at the Institut Mittag-Leffler in Djursholm, Sweden. We are grateful to the institute and the organisers of the workshop for the wonderful working environment.

\bibliography{starbib}

@article {Tu61,
    AUTHOR = {Tutte, W. T.},
     TITLE = {On the problem of decomposing a graph into {$n$} connected
              factors},
   JOURNAL = {J. London Math. Soc.},
  FJOURNAL = {The Journal of the London Mathematical Society},
    VOLUME = {36},
      YEAR = {1961},
     PAGES = {221--230},
      ISSN = {0024-6107,1469-7750},
   MRCLASS = {05.45},
  MRNUMBER = {140438},
MRREVIEWER = {G.\ A.\ Dirac},
       DOI = {10.1112/jlms/s1-36.1.221},
       URL = {https://doi.org/10.1112/jlms/s1-36.1.221},
}

@article {NW61,
    AUTHOR = {Nash-Williams, C. St.\ J. A.},
     TITLE = {Edge-disjoint spanning trees of finite graphs},
   JOURNAL = {J. London Math. Soc.},
  FJOURNAL = {The Journal of the London Mathematical Society},
    VOLUME = {36},
      YEAR = {1961},
     PAGES = {445--450},
      ISSN = {0024-6107,1469-7750},
   MRCLASS = {05.45},
  MRNUMBER = {133253},
MRREVIEWER = {W.\ T.\ Tutte},
       DOI = {10.1112/jlms/s1-36.1.445},
       URL = {https://doi.org/10.1112/jlms/s1-36.1.445},
}

@article {SchiermeyerSotak,
    AUTHOR = {Schiermeyer, Ingo and Sot\'ak, Roman},
     TITLE = {Rainbow numbers for graphs containing small cycles},
   JOURNAL = {Graphs Combin.},
  FJOURNAL = {Graphs and Combinatorics},
    VOLUME = {31},
      YEAR = {2015},
    NUMBER = {6},
     PAGES = {1985--1991},
      ISSN = {0911-0119,1435-5914},
   MRCLASS = {05C55 (05C15 05C35)},
  MRNUMBER = {3417209},
MRREVIEWER = {Colton\ Magnant},
       DOI = {10.1007/s00373-015-1577-7},
       URL = {https://doi.org/10.1007/s00373-015-1577-7},
}

@article{ErdosRado,
  title={A combinatorial theorem},
  author={Erd{\H{o}}s, Paul and Rado, Richard},
JOURNAL = {J. London Math. Soc.},
  volume={1},
  number={4},
  pages={249--255},
  year={1950},
  publisher={Wiley Online Library}
}

@article{ErdosStone,
    AUTHOR = {Erd\H{o}s, P. and Stone, A. H.},
     TITLE = {On the structure of linear graphs},
   JOURNAL = {Bull. Amer. Math. Soc.},
  FJOURNAL = {Bulletin of the American Mathematical Society},
    VOLUME = {52},
      YEAR = {1946},
     PAGES = {1087--1091},
      ISSN = {0002-9904},
   MRCLASS = {56.0X},
       DOI = {10.1090/S0002-9904-1946-08715-7},
       URL = {https://doi.org/10.1090/S0002-9904-1946-08715-7},
}

@article{RainbowSurvey,
    AUTHOR = {Fujita, Shinya and Magnant, Colton and Ozeki, Kenta},
     TITLE = {Rainbow generalizations of {R}amsey theory: a survey},
   JOURNAL = {Graphs Combin.},
  FJOURNAL = {Graphs and Combinatorics},
    VOLUME = {26},
      YEAR = {2010},
    NUMBER = {1},
     PAGES = {1--30},
      ISSN = {0911-0119,1435-5914},
   MRCLASS = {05C55 (05C15 05C35)},
       DOI = {10.1007/s00373-010-0891-3},
       URL = {https://doi.org/10.1007/s00373-010-0891-3},
}

@incollection{PokrovskiySurvey,
    AUTHOR = {Pokrovskiy, Alexey},
     TITLE = {Rainbow subgraphs and their applications},
 BOOKTITLE = {Surveys in combinatorics 2022},
    SERIES = {London Math. Soc. Lecture Note Ser.},
    VOLUME = {481},
     PAGES = {191--214},
 PUBLISHER = {Cambridge Univ. Press, Cambridge},
      YEAR = {2022},
      ISBN = {978-1-009-09622-5},
   MRCLASS = {05C55 (05C15)},
  MRNUMBER = {4421403},
}

@article{AlonPokrovskiySudakov,
    AUTHOR = {Alon, Noga and Pokrovskiy, Alexey and Sudakov, Benny},
     TITLE = {Random subgraphs of properly edge-coloured complete graphs and
              long rainbow cycles},
   JOURNAL = {Israel J. Math.},
  FJOURNAL = {Israel Journal of Mathematics},
    VOLUME = {222},
      YEAR = {2017},
    NUMBER = {1},
     PAGES = {317--331},
      ISSN = {0021-2172,1565-8511},
   MRCLASS = {05C15 (05C38 05C45 05C80)},
       DOI = {10.1007/s11856-017-1592-x},
       URL = {https://doi.org/10.1007/s11856-017-1592-x},
}

@article{MontgomeryPokrovskiySudakov,
    AUTHOR = {Montgomery, Richard and Pokrovskiy, Alexey and Sudakov,
              Benjamin},
     TITLE = {Decompositions into spanning rainbow structures},
   JOURNAL = {Proc. London Math. Soc. (3)},
  FJOURNAL = {Proceedings of the London Mathematical Society. Third Series},
    VOLUME = {119},
      YEAR = {2019},
    NUMBER = {4},
     PAGES = {899--959},
      ISSN = {0024-6115,1460-244X},
   MRCLASS = {05B15},
       DOI = {10.1112/plms.12245},
       URL = {https://doi.org/10.1112/plms.12245},
}

@incollection{ErdosSimonovitsSos,
    AUTHOR = {Erd\H{o}s, P. and Simonovits, M. and S\'os, V. T.},
     TITLE = {Anti-{R}amsey theorems},
 BOOKTITLE = {Infinite and finite sets ({C}olloq., {K}eszthely, 1973;
              dedicated to {P}. {E}rd\H os on his 60th birthday), {V}ols.
              {I}, {II}, {III}},
    SERIES = {Colloq. Math. Soc. J\'anos Bolyai},
    VOLUME = {Vol. 10},
     PAGES = {633--643},
 PUBLISHER = {North-Holland, Amsterdam-London},
      YEAR = {1975},
   MRCLASS = {05C15},
}

@article{MontellanoNeumann,
    AUTHOR = {Montellano-Ballesteros, J. J. and Neumann-Lara, V.},
     TITLE = {An anti-{R}amsey theorem},
   JOURNAL = {Combinatorica},
  FJOURNAL = {Combinatorica. An International Journal on Combinatorics and
              the Theory of Computing},
    VOLUME = {22},
      YEAR = {2002},
    NUMBER = {3},
     PAGES = {445--449},
      ISSN = {0209-9683,1439-6912},
   MRCLASS = {05C55 (05C35)},
       DOI = {10.1007/s004930200023},
       URL = {https://doi.org/10.1007/s004930200023},
}

@article{AxenovichIverson,
    AUTHOR = {Axenovich, Maria and Iverson, Perry},
     TITLE = {Edge-colorings avoiding rainbow and monochromatic subgraphs},
   JOURNAL = {Discrete Math.},
  FJOURNAL = {Discrete Mathematics},
    VOLUME = {308},
      YEAR = {2008},
    NUMBER = {20},
     PAGES = {4710--4723},
      ISSN = {0012-365X,1872-681X},
   MRCLASS = {05C55 (05C15)},
       DOI = {10.1016/j.disc.2007.08.092},
       URL = {https://doi.org/10.1016/j.disc.2007.08.092},
}

@article{JamisonJiangLing,
    AUTHOR = {Jamison, Robert E. and Jiang, Tao and Ling, Alan C. H.},
     TITLE = {Constrained {R}amsey numbers of graphs},
   JOURNAL = {J. Graph Theory},
  FJOURNAL = {Journal of Graph Theory},
    VOLUME = {42},
      YEAR = {2003},
    NUMBER = {1},
     PAGES = {1--16},
      ISSN = {0364-9024,1097-0118},
   MRCLASS = {05C55},
       DOI = {10.1002/jgt.10072},
       URL = {https://doi.org/10.1002/jgt.10072},
}

@article{MaamounMeyniel,
    AUTHOR = {Maamoun, M. and Meyniel, H.},
     TITLE = {On a problem of {G}. {H}ahn about coloured {H}amiltonian paths
              in {$K\sb{2n}$}},
   JOURNAL = {Discrete Math.},
  FJOURNAL = {Discrete Mathematics},
    VOLUME = {51},
      YEAR = {1984},
    NUMBER = {2},
     PAGES = {213--214},
      ISSN = {0012-365X,1872-681X},
   MRCLASS = {05C15 (05C45)},
       DOI = {10.1016/0012-365X(84)90073-6},
       URL = {https://doi.org/10.1016/0012-365X(84)90073-6},
}

@article{Burr,
    AUTHOR = {Burr, Stefan A.},
     TITLE = {An inequality involving the vertex arboricity and edge
              arboricity of a graph},
   JOURNAL = {J. Graph Theory},
  FJOURNAL = {Journal of Graph Theory},
    VOLUME = {10},
      YEAR = {1986},
    NUMBER = {3},
     PAGES = {403--404},
      ISSN = {0364-9024,1097-0118},
   MRCLASS = {05C15 (05C99)},
       DOI = {10.1002/jgt.3190100315},
       URL = {https://doi.org/10.1002/jgt.3190100315},
}

@article{Redei,
  title={Ein kombinatorischer {S}atz},
  author={R{\'e}dei, L{\'a}szl{\'o}},
  journal={Acta Litt. Szeged},
  volume={7},
  number={39-43},
  pages={97},
  year={1934}
}

@article{EhardGlockJoos,
    AUTHOR = {Ehard, Stefan and Glock, Stefan and Joos, Felix},
     TITLE = {A rainbow blow-up lemma for almost optimally bounded
              edge-colourings},
   JOURNAL = {Forum Math. Sigma},
  FJOURNAL = {Forum of Mathematics. Sigma},
    VOLUME = {8},
      YEAR = {2020},
     PAGES = {Paper No. e37, 32},
      ISSN = {2050-5094},
   MRCLASS = {05C15 (05B40 05C35 05C51 05C60 05C70 05C78)},
       DOI = {10.1017/fms.2020.38},
       URL = {https://doi.org/10.1017/fms.2020.38},
}

@article{Andersen,
    AUTHOR = {Andersen, Lars D\o vling},
     TITLE = {Hamilton circuits with many colours in properly edge-coloured
              complete graphs},
   JOURNAL = {Math. Scand.},
  FJOURNAL = {Mathematica Scandinavica},
    VOLUME = {64},
      YEAR = {1989},
    NUMBER = {1},
     PAGES = {5--14},
      ISSN = {0025-5521,1903-1807},
   MRCLASS = {05C15 (05C45)},
       DOI = {10.7146/math.scand.a-12245},
       URL = {https://doi.org/10.7146/math.scand.a-12245},
}

@article{BukhConlon,
    AUTHOR = {Bukh, Boris and Conlon, David},
     TITLE = {Rational exponents in extremal graph theory},
   JOURNAL = {J. Eur. Math. Soc. (JEMS)},
  FJOURNAL = {Journal of the European Mathematical Society (JEMS)},
    VOLUME = {20},
      YEAR = {2018},
    NUMBER = {7},
     PAGES = {1747--1757},
      ISSN = {1435-9855,1435-9863},
   MRCLASS = {05C35},
       DOI = {10.4171/JEMS/798},
       URL = {https://doi.org/10.4171/JEMS/798},
}

@article{ConlonJanzer,
    AUTHOR = {Conlon, David and Janzer, Oliver},
     TITLE = {Rational exponents near two},
   JOURNAL = {Adv. Comb.},
  FJOURNAL = {Advances in Combinatorics},
      YEAR = {2022},
     PAGES = {Paper No. 9, 10},
      ISSN = {2517-5599},
   MRCLASS = {05C35},
  MRNUMBER = {4527790},
}

@article{JiangQiu,
    AUTHOR = {Jiang, Tao and Qiu, Yu},
     TITLE = {Many {T}ur\'an exponents via subdivisions},
   JOURNAL = {Combin. Probab. Comput.},
  FJOURNAL = {Combinatorics, Probability and Computing},
    VOLUME = {32},
      YEAR = {2023},
    NUMBER = {1},
     PAGES = {134--150},
      ISSN = {0963-5483,1469-2163},
   MRCLASS = {05C35},
       DOI = {10.1017/s0963548322000177},
       URL = {https://doi.org/10.1017/s0963548322000177},
}

@incollection{FurediSimonovits,
    AUTHOR = {F\"uredi, Zolt\'an and Simonovits, Mikl\'os},
     TITLE = {The history of degenerate (bipartite) extremal graph problems},
 BOOKTITLE = {Erd\H{o}s centennial},
    SERIES = {Bolyai Soc. Math. Stud.},
    VOLUME = {25},
     PAGES = {169--264},
 PUBLISHER = {J\'anos Bolyai Math. Soc., Budapest},
      YEAR = {2013},
      ISBN = {978-963-9453-18-0; 978-3-642-39285-6},
   MRCLASS = {05-02 (01A70)},
       DOI = {10.1007/978-3-642-39286-3\_7},
       URL = {https://doi.org/10.1007/978-3-642-39286-3_7},
}

@article{KovariSosTuran,
    AUTHOR = {K\H{o}vari, T. and S\'os, V. T. and Tur\'an, P.},
     TITLE = {On a problem of {K}. {Z}arankiewicz},
   JOURNAL = {Colloq. Math.},
  FJOURNAL = {Colloquium Mathematicum},
    VOLUME = {3},
      YEAR = {1954},
     PAGES = {50--57},
      ISSN = {0010-1354,1730-6302},
   MRCLASS = {27.2X},
       DOI = {10.4064/cm-3-1-50-57},
       URL = {https://doi.org/10.4064/cm-3-1-50-57},
}

@article{KollarRonyaiSzabo,
    AUTHOR = {Koll\'ar, J\'anos and R\'onyai, Lajos and Szab\'o, Tibor},
     TITLE = {Norm-graphs and bipartite {T}ur\'an numbers},
   JOURNAL = {Combinatorica},
  FJOURNAL = {Combinatorica. An International Journal on Combinatorics and
              the Theory of Computing},
    VOLUME = {16},
      YEAR = {1996},
    NUMBER = {3},
     PAGES = {399--406},
      ISSN = {0209-9683},
   MRCLASS = {05C35},
  MRNUMBER = {1417348},
MRREVIEWER = {W.\ G.\ Brown},
       DOI = {10.1007/BF01261323},
       URL = {https://doi.org/10.1007/BF01261323},
}

@article{Moon,
    AUTHOR = {Moon, J. W.},
     TITLE = {On subtournaments of a tournament},
   JOURNAL = {Canad. Math. Bull.},
  FJOURNAL = {Canadian Mathematical Bulletin. Bulletin Canadien de
              Math\'ematiques},
    VOLUME = {9},
      YEAR = {1966},
     PAGES = {297--301},
      ISSN = {0008-4395,1496-4287},
   MRCLASS = {05.60},
       DOI = {10.4153/CMB-1966-038-7},
       URL = {https://doi.org/10.4153/CMB-1966-038-7},
}

@article{AxenovichKundgen,
    AUTHOR = {Axenovich, Maria and K\"undgen, Andr\'e},
     TITLE = {On a generalized anti-{R}amsey problem},
   JOURNAL = {Combinatorica},
  FJOURNAL = {Combinatorica. An International Journal on Combinatorics and
              the Theory of Computing},
    VOLUME = {21},
      YEAR = {2001},
    NUMBER = {3},
     PAGES = {335--349},
      ISSN = {0209-9683,1439-6912},
   MRCLASS = {05C55 (05C35 05D40)},
       DOI = {10.1007/s004930100000},
       URL = {https://doi.org/10.1007/s004930100000},
}

@article{JungicKaiserKral,
    AUTHOR = {Jungi\'c, Veselin and Kaiser, Tom\'a\us and Kr\'al, Daniel},
     TITLE = {A note on edge-colourings avoiding rainbow {$K_4$} and
              monochromatic {$K_m$}},
   JOURNAL = {Electron. J. Combin.},
  FJOURNAL = {Electronic Journal of Combinatorics},
    VOLUME = {16},
      YEAR = {2009},
    NUMBER = {1},
     PAGES = {Note 19, 9},
      ISSN = {1077-8926},
       DOI = {10.37236/257},
       URL = {https://doi.org/10.37236/257},
}

@article{Li,
    AUTHOR = {Li, Hao},
     TITLE = {Rainbow {$C_3$}'s and {$C_4$}'s in edge-colored graphs},
   JOURNAL = {Discrete Math.},
  FJOURNAL = {Discrete Mathematics},
    VOLUME = {313},
      YEAR = {2013},
    NUMBER = {19},
     PAGES = {1893--1896},
      ISSN = {0012-365X,1872-681X},
   MRCLASS = {05C15 (05C38)},
       DOI = {10.1016/j.disc.2012.11.024},
       URL = {https://doi.org/10.1016/j.disc.2012.11.024},
}

@article{CzygrinowMollaNagle,
  title={Rainbow cliques in edge-colored graphs},
  author={Czygrinow, Andrzej and Molla, Theodore and Nagle, Brendan},
  journal={arXiv preprint arXiv:2407.08098},
  year={2024}
}

@article{FoxSudakov,
    AUTHOR = {Fox, Jacob and Sudakov, Benny},
     TITLE = {Dependent random choice},
   JOURNAL = {Random Struct. Alg.},
  FJOURNAL = {Random Structures \& Algorithms},
    VOLUME = {38},
      YEAR = {2011},
    NUMBER = {1-2},
     PAGES = {68--99},
      ISSN = {1042-9832,1098-2418},
   MRCLASS = {05D40 (05C35 05C55 05D10 60C05)},
  MRNUMBER = {2768884},
MRREVIEWER = {Hamed\ Hatami},
       DOI = {10.1002/rsa.20344},
       URL = {https://doi.org/10.1002/rsa.20344},
}

@article{Furedi,
    AUTHOR = {F\"uredi, Zolt\'an},
     TITLE = {New asymptotics for bipartite {T}ur\'an numbers},
   JOURNAL = {J. Combin. Theory Ser. A},
  FJOURNAL = {Journal of Combinatorial Theory. Series A},
    VOLUME = {75},
      YEAR = {1996},
    NUMBER = {1},
     PAGES = {141--144},
      ISSN = {0097-3165,1096-0899},
   MRCLASS = {05C35 (05C15 05C50)},
       DOI = {10.1006/jcta.1996.0067},
       URL = {https://doi.org/10.1006/jcta.1996.0067},
}

@article{AlonHooryLinial,
     AUTHOR = {Alon, Noga and Hoory, Shlomo and Linial, Nathan},
     TITLE = {The {M}oore bound for irregular graphs},
   JOURNAL = {Graphs Combin.},
  FJOURNAL = {Graphs and Combinatorics},
    VOLUME = {18},
      YEAR = {2002},
    NUMBER = {1},
     PAGES = {53--57},
      ISSN = {0911-0119,1435-5914},
   MRCLASS = {05C35},
       DOI = {10.1007/s003730200002},
       URL = {https://doi.org/10.1007/s003730200002},
}

@article{Benson,
    AUTHOR = {Benson, Clark T.},
     TITLE = {Minimal regular graphs of girths eight and twelve},
   JOURNAL = {Canadian J. Math.},
  FJOURNAL = {Canadian Journal of Mathematics. Journal Canadien de
              Math\'ematiques},
    VOLUME = {18},
      YEAR = {1966},
     PAGES = {1091--1094},
      ISSN = {0008-414X,1496-4279},
   MRCLASS = {05.40},
       DOI = {10.4153/CJM-1966-109-8},
       URL = {https://doi.org/10.4153/CJM-1966-109-8},
}

@article{Singleton,
    AUTHOR = {Singleton, Robert},
     TITLE = {On minimal graphs of maximum even girth},
   JOURNAL = {J. Combin. Theory},
  FJOURNAL = {Journal of Combinatorial Theory},
    VOLUME = {1},
      YEAR = {1966},
     PAGES = {306--332},
      ISSN = {0021-9800},
   MRCLASS = {05.40},
}

@article{Reiman,
    AUTHOR = {Reiman, I.},
     TITLE = {{\"U}ber ein {P}roblem von {K}. {Z}arankiewicz},
   JOURNAL = {Acta Math. Acad. Sci. Hungar.},
  FJOURNAL = {Acta Mathematica. Academiae Scientiarum Hungaricae},
    VOLUME = {9},
      YEAR = {1958},
     PAGES = {269--273},
      ISSN = {0001-5954,1588-2632},
   MRCLASS = {05.00},
       DOI = {10.1007/BF02020254},
       URL = {https://doi.org/10.1007/BF02020254},
}

@book{janson2011random,
  title={Random graphs},
  author={Janson, Svante and {\L}uczak, Tomasz and Ruci\'nski, Andrzej},
  year={2011},
  publisher={John Wiley \& Sons}
}
\bibliographystyle{abbrv}

\end{document}